\documentclass[12pt,tbtags,leqno]{amsart}

\usepackage{amssymb}
\usepackage{amsthm}
\usepackage{amsmath}
\usepackage{verbatim}
\usepackage[numbers,sort&compress]{natbib}

\numberwithin{equation}{section}

\newcommand{\HH}{\mathcal{H}}
\newcommand{\HK}{\mathcal{K}}
\newcommand{\HT}{\mathcal{T}}
\newcommand{\HM}{\mathcal{M}}

\newcommand{\HG}{\mathcal{G}}

\newcommand{\HS}{\mathcal{S}}

\newcommand{\HD}{\mathcal{D}}
\newcommand{\HB}{\mathcal{B}}
\newcommand{\HN}{\mathcal{N}}

\newcommand{\D}{\mathbb{D}}

\newcommand{\C}{\mathbb{C}}

\newcommand{\ran}{\mathrm{ran \ }}

\newcommand{\rank}{\mathrm{rank \ }}

\newcommand{\Concave}{\mathfrak{C}}
\newcommand{\ConcaveNorm}{\mathfrak{C}_{\mathrm{norm}}}
\newcommand{\ConcaveNil}{\mathfrak{C}_{\mathrm{nil}}}

\theoremstyle{plain}
\newtheorem{theorem}{Theorem}[section]
\newtheorem{lemma}[theorem]{Lemma}
\newtheorem{remark}[theorem]{Remark}

\newtheorem{proposition}[theorem]{Proposition}
\newtheorem{corollary}[theorem]{Corollary}

\newtheorem{definition}[theorem]{Definition}
\newtheorem{question}[theorem]{Question}
\theoremstyle{definition}

\makeatletter
\@namedef{subjclassname@2020}{%
  \textup{2020} Mathematics Subject Classification}
\makeatother

\begin{document}
\title{Concave shifts and the complete Nevanlinna--Pick property}

\author{Shuaibing Luo and Eskil Rydhe}

\address{School of Mathematics, Hunan University, Changsha, Hunan, 410082, PR China}
\email{sluo@hnu.edu.cn}

\address{Centre for Mathematical Sciences, Lund University, Box 118, SE-221 00 Lund, Sweden}
\email{eskil.rydhe@math.lu.se}



\date{}

\maketitle

\begin{abstract}
Let $\HH$ be a non-zero reproducing kernel Hilbert space of holomorphic functions on the unit disc, and assume $\HH$ is invariant under the forward and backward shifts. We show that if the forward shift on $\HH$ is concave, then $\HH$ is a complete Nevanlinna--Pick (CNP) space. This answers a question by Shimorin, and greatly extends his result that any space of Dirichlet type is CNP. As a byproduct of our investigation, we also show that if the forward shift on $\HH$ is concave, then any multiplier of $\HH$ is the SOT limit of its Cesàro means. As applications, we prove a commutant lifting theorem for concave cyclic analytic operators, and several results about shift-invariant subspaces in de~Branges--Rovnyak spaces.
\end{abstract}

\section{Introduction}

Let $\HH$ be a complex separable Hilbert space, with $\HB(\HH)$ the corresponding space of bounded linear operators. Unless we explicitly state otherwise, $\HH$ will be a non-zero \emph{reproducing kernel Hilbert space} (RKHS) of holomorphic functions $f\colon\D\to\C$. A brief background on RKHSs is provided in Section~\ref{sec:Preliminaries}. We generally assume that \(\mathcal H\) is invariant under the forward and backward shifts, respectively defined by
\[
\left(M_zf\right)(z)=zf(z),\qquad \left(Lf\right)(z)=\frac{f(z)-f(0)}{z}.
\]
We write \(T= (M_z, \HH)\) for the operator on $\HH$ induced by $M_z$. By the closed graph theorem, $T$ is bounded. Similarly, $(L,\HH)$ is also bounded.

The \emph{reproducing kernel} of $\HH$ will be denoted $K$, or sometimes $K^\HH$. If there exists $\alpha\in\D$ such that $K(\alpha,\alpha)\ne 0$, and
\begin{align}\label{eq:CNP}
	F(z,w)=1-\frac{K(z,\alpha)K(\alpha,w)}{K(\alpha,\alpha)K(z,w)}
\end{align}
defines a positive semi-definite function on $\D\times\D$, then we say that $K$ is a \emph{complete Nevanlinna--Pick} (CNP) kernel. If $K$ is a CNP kernel, then we say that $\HH$ is CNP. The notion of CNP spaces originates from the classical Pick interpolation problem in the Hardy space $H^2$, \cite{Pick1915}. In the so-called operator-theoretic approach, pioneered by Sarason~\cite{Sarason1967}, interpolation problems are reformulated in terms of existence of multipliers. Sarason's approach quickly led Sz.-Nagy and Foias to formulate their commutant lifting theorem for contractions on Hilbert spaces~\cite{SzNagyFoias1968a,SzNagyFoias1968b}, effectively initiating a new subfield of operator theory. Other developments revealed that the CNP property in itself gives rise to a rich function theory. As a striking example of this, we mention a result by Aleman--Hartz--McCarthy--Richter~\cite{AHMR17}, that if $\HH$ is a normalized CNP space, then any function in $\HH$ is in fact the quotient of two multipliers on $\HH$. For $H^2$, the prototypical example of a CNP space, this is just the classical theorem that any function in $H^2$ is a quotient of two functions in $H^\infty$. A systematic treatment of Pick interpolation, including many historical notes, is given in \cite{AM02}. The connection between Pick interpolation and \eqref{eq:CNP} is known as the McCullough--Quiggin theorem~\cite{McCullough1994,Qui93}, see also~\cite{AM00} and~\cite[Theorem~7.6]{AM02}. We remark that the technical condition of $K$ being ``irreducible'' is always satisfied in the examples we consider. See Section~\ref{sec:Preliminaries}, and the discussion following Lemma~\ref{hbmodsscnprhB}.

Although the function theory for CNP spaces is by now quite developed, e.g.~\cite{AHMR17,AHMR18,AHMR19,AHMR21,AHMR23}, it is often challenging to verify whether a given space $\HH$ is CNP. We illustrate this with two common examples.

The so-called \emph{de~Branges--Rovnyak spaces} $\HH[B]$ introduced in \cite{debranges} continue to generate attention, see for instance~\cite{AlemanMalman,BallBoloBasics,CGL24,FM162}. The spaces themselves are rather abstract in nature, but the kernels are explicitly known. Nevertheless, the direct analysis of \eqref{eq:CNP} to determine whether $\HH[B]$ is CNP is non-trivial even for the so-called rank-one case, see~\cite{Chu20}. For a general $\HH[B]$ of arbitrary or infinite rank, it appears to be a very difficult problem to verify the CNP property.

The so-called \emph{Dirichlet type spaces} $D(\mu)$ were introduced by Aleman~\cite{Al93}, following important contributions by Richter~\cite{Ri88,Ri91}. These are concrete spaces, but the kernels are implicitly defined by the Riesz representation theorem, and apart from a handful of examples the kernels have no explicit representations. This necessitates using rather indirect methods for studying \eqref{eq:CNP}.

In 2002, Shimorin proved the remarkable theorem that every $D(\mu)$ space is in fact CNP~\cite[Theorem 1.1]{Sh02}. The proof does not shed any light on what the reproducing kernels are. In the same paper, Shimorin suggested that a much larger class of RKHSs could be CNP. To make this precise, consider for a moment a general Hilbert space $\HH$. $T\in\HB(\HH)$ is called \emph{concave} if
\begin{align}\label{tstarrel3}
	T^{*2}T^2 - 2T^*T + I \leq 0,
\end{align}
and \emph{logarithmically concave} if
\begin{align}\label{tstarrel}
	T^{*2}T^2 -T^*TT^*T \le 0.
\end{align}
Introducing the standard notation \(\Delta=T^*T-I\), the above conditions are equivalent to \(T^*\Delta T\leq\Delta\), and \(T^*\Delta T\leq\Delta + \Delta^2\), respectively. In particular, every concave operator is logarithmically concave.  It is known that concave operators are \emph{expansive}, meaning $\Delta \ge 0$, e.g.~\cite[Lemma~1]{Ri88}. Any operator of the form $T=(M_z,D(\mu))$ is concave. In Section~\ref{sec:Comparison}, we present a few illustrations of why $D(\mu)$ is a very particular example of a RKHS with concave forward shift.
\begin{question}[Shimorin~\cite{Sh02},~Section~5]\label{Shimorin}
	Let $\HH$ be a non-zero RKHS of holomorphic functions $f\colon\D\to\C$. Assume that $\HH$ is invariant under $M_z$ and $L$. For $T=(M_z,\HH)$, does either one of the conditions (\ref{tstarrel3}) or  (\ref{tstarrel}) imply that $\HH$ is a CNP space?
\end{question}

Our main result is a complete resolution of Shimorin's question. The answer is negative for logarithmically concave operators, and affirmative for concave operators.
\begin{theorem}\label{mainthmonelogconnot}
	There exists $\HH$, a RKHS of holomorphic functions $f\colon\D\to\C$, with the following properties.
	\begin{enumerate}
		\item $\HH$ is invariant under $M_z$ and $L$.
		\item $T=(M_z,\HH)$ is expansive and logarithmically concave.
		\item $\Delta = (T^*T-I)$ has rank $3$.
		\item $\HH$ is not CNP.
	\end{enumerate}
\end{theorem}
\begin{theorem}\label{mainthmconcave}
	Let $\HH$ be a RKHS of holomorphic functions $f\colon\D\to\C$. Assume that $\HH$ is invariant under $M_z$ and $L$. If $T = (M_z, \HH)$ is concave, then $\HH$ is CNP.
\end{theorem}

Theorem~\ref{mainthmonelogconnot} is proved by construction of an explicit example. Before we indicate the proof of Theorem~\ref{mainthmconcave}, we introduce some terminology. We write $P_{\HM}$ for the orthogonal projection onto a closed subspace $\HM\subseteq\HH$. We write $[E]_T$ to indicate the smallest closed $T$-invariant subspace containing a set $E$.
\begin{definition}
	Let $\HH$ be an RKHS of holomorphic functions $f\colon\D\to\C$. We say that the reproducing kernel $K^\HH$ is \emph{normalized} if $K^\HH(z,0)=1$ for every $z\in\D$. Assume $\HH$ is invariant under $M_z$ and $L$. Let $T=(M_z,\HH)$, $\Delta = T^*T-I$, $\HN = [\ran \Delta]_{T^*}$, and $A=P_\HN T|_\HN$.
	\begin{enumerate}
		\item We write $\HH\in\Concave$ to indicate that $T$ is concave.
		\item We write $\HH\in\ConcaveNorm$ to indicate that $\HH\in\Concave$, and $K^\HH$ is normalized.
		\item We write $\HH\in\ConcaveNil$ to indicate that $\HH\in\mathfrak{C}_{\mathrm{norm}}$, $\dim\HN<\infty$, and $A$ is nilpotent.
	\end{enumerate}
\end{definition}
\begin{remark}\label{remark:NisranDelta}
	If $T$ is logarithmically concave, then $T^*\Delta T\le \Delta + \Delta^2$, so if $T$ is also expansive, then $\ker\Delta$ is $T$-invariant, and $\overline{\ran\Delta} = (\ker \Delta)^\perp$ is $T^*$-invariant. In particular, $\HH\in\mathfrak{C}_{\mathrm{nil}}$ implies that $\HN=\ran\Delta$.
\end{remark}

The hypothesis of Theorem~\ref{mainthmconcave} is that $\HH\in\mathfrak{C}$. Since concave operators are expansive, \cite[Theorem~4.6]{LGR} by Gu, Richter, and the first author allows for a non-trivial but straightforward reduction to the case $\HH\in\mathfrak{C}_{\mathrm{norm}}$. A more involved part of the present paper is a further reduction to the hypothesis $\HH\in\mathfrak{C}_{\mathrm{nil}}$. Proposition~\ref{nilpotconcnp} asserts that $\HH$ is CNP whenever $\HH\in\mathfrak{C}_{\mathrm{nil}}$. The main idea of the proof is to identify a shift-invariant subspace $\HM\subseteq\HH$ on which $(M_z,\HM)$ unitarily equivalent to $(M_z,H^2)$, and then apply a finite rank extension scheme which preserves the CNP property. Combining Proposition~\ref{nilpotconcnp} with an approximation argument yields that $\HH$ is CNP whenever $\HH\in\mathfrak{C}_\mathrm{norm}$. Lemma~\ref{hbmodsscnprhB} then implies Theorem~\ref{mainthmconcave}.

A byproduct of our investigation is what appears to be a novel description of certain de~Branges--Rovnyak spaces. Let $E$ be an auxiliary Hilbert space, $C\in\HB(E)$ a contraction, and $v\in E$. For a given $f\colon z\mapsto \sum_{n\ge 0}f_nz^n\in H^2$, consider the equations
\begin{equation}\label{fysaiedn}
	y_k - C y_{k+1} = f_k v, \qquad k = 1, 2, 3, \ldots,
\end{equation}
with unknowns $y_1, y_2, \ldots \in E$. If a solution $y = (y_1, y_2, \ldots) \in \ell^2(\mathbb{N}, E)$ exists, then it is unique. Indeed, if $z$ is the difference of two solutions, then $z_k - C z_{k+1} = 0$ for $k \geq 1$. It follows that $z_k = C^m z_{k+m}$ for every $m \ge 1$. Since $\|C\| \le 1$ and $z_{k+m} \to 0$, this implies $z_k = 0$ for every $k$.
\begin{definition}\label{def:G}
	With notation as in the previous paragraph, the Hilbert space $\mathcal{G}=\mathcal{G}(E,C,v)$ consists of all $f\in H^2$ for which \eqref{fysaiedn} has a solution $y=y(f)\in\ell^2(\mathbb{N}, E)$. The norm on $\mathcal{G}$ is defined by
	\begin{equation*}
		\|f\|_{\mathcal{G}}^2 = \|f\|_{H^2}^2 + \|y(f)\|_{\ell^2(E)}^2.
	\end{equation*}
\end{definition}

\begin{proposition}\label{prop:Gspaces}
	If $\HH[B]\in\ConcaveNorm$, then $\HH[B]=\HG(E,C,v)$, where $E=\overline{\ran \Delta^{1/2}}$, $C\in\HB(E)$ is densely defined by $C\Delta^{1/2}f=\Delta^{1/2} T f$, and $v = \Delta^{1/2} 1$.
\end{proposition}

The above description of $\HH[B]\in\ConcaveNorm$ is not only useful for the proof of Theorem~\ref{mainthmconcave}, but also allows us to deduce new results about multipliers on $\HH[B]$. For any holomorphic functions $f,\psi\colon\D\to\C$, define $M_\psi f := \psi f$, and recall that the multiplier algebra of a RKHS $\HH$ is defined by
\begin{align}\label{eq:MultiplierAlgebra}
	\operatorname{Mult}(\mathcal{H}) = \{\psi \in \text{Hol}(\D) \mid M_\psi\HH\subseteq\HH\}.
\end{align}
If $\psi\in\operatorname{Mult}(\HH)$, then $M_\psi$ is bounded according to the closed graph theorem. The norm of $\psi\in\operatorname{Mult}(\HH)$ is the operator norm of $(M_\psi,\HH)$. Proposition~\ref{prop:Gspaces} allows us to prove that if $\HH\in\Concave$, then any $\psi\in\operatorname{Mult}(\HH)$ is the SOT-limit of its Cesàro means. In particular, polynomials are SOT-dense in the multiplier algebra. Recall that if $\psi\colon z\mapsto \sum_{n\ge 0}\psi_nz^n$, then its Cesàro means are given by
\begin{align*}
	(\sigma_n \psi) (z)=\sum_{k=0}^n\left(1-\frac{k}{n+1}\right)\psi_k z^k.
\end{align*}
\begin{theorem}\label{thm:Cesaro}
	If $\HH\in\Concave$ and $\psi\in\operatorname{Mult}(\HH)$, then $M_{\sigma_n\psi}\to M_\psi$ with convergence in the strong operator topology.
\end{theorem}
In certain applications, e.g. \cite[Theorem~1.5]{Har23}, \cite[Theorem~1.1]{Sh05}, hypotheses are formulated in terms of subspaces being invariant for all multipliers. Theorem~\ref{thm:Cesaro} implies that if $\HH\in\Concave$, then it suffices to consider the single multiplication operators $M_z$. Imposing a circular symmetry assumption on the kernel, this is true in much greater generality~\cite[Proposition~2.5]{CHS2020}. Without the symmetry assumption, \cite[Lemma~5.4]{RichterSundberg1991} readily implies that polynomials are WOT-dense in $\operatorname{Mult}(D(\mu))$ provided that $\mu$ is supported on $\partial\D$. Theorem~\ref{thm:Cesaro} extends this result to a much larger class of function spaces.

The paper is organized as follows. Section~\ref{sec:Preliminaries} contains some preliminaries on RKHSs, in particular $\HH[B]$ and $D(\mu)$. In Section~\ref{sec:Comparison}, we relate Dirichlet type spaces to the condition $\HH\in\Concave$, and give a few illustrations of why the latter condition is more general. In Section~\ref{sec:ProofMainThm2}, we prove Theorem \ref{mainthmonelogconnot} by constructing an explicit RKHS which is not CNP but still has logarithmically concave shift. In Section~\ref{sec:ProofMainThm} we prove Theorem~\ref{mainthmconcave} and the auxiliary result Proposition~\ref{prop:Gspaces}. Theorem~\ref{thm:Cesaro} is proved in Section~\ref{sec:Multipliers}. In Section~\ref{sec:Applications} we demonstrate some applications,  including a commutant lifting theorem for concave cyclic analytic operators, and structure theorems for $M_z$-invariant subspaces of $\HH\in\Concave$. 

\section{Preliminaries}\label{sec:Preliminaries}

As in the introduction, let $\HH$ be a Hilbert space of holomorphic functions $f\colon \D\to \C$. For a given $w\in\D$, if the evaluation $f\mapsto f(w)$ defines a bounded linear functional on $\HH$, then there exists $k_w\in\HH$ such that $f(w)=\langle f,k_w\rangle_\HH$ whenever $f\in\HH$. If $k_w$ exists for every $w\in\D$, then we call $\HH$ a \emph{reproducing kernel Hilbert space} (RKHS). If $\HH$ is a RKHS, then the sesquiholomorphic function $K\colon\D\times\D\to\C$ defined by $K(z,w):=k_w(z)$ is called the \emph{reproducing kernel} of $\HH$. A good general reference on RKHS is~\cite{PaulsenRag}.

The space $\HH$ is called \emph{irreducible} if for any two distinct $w_1,w_2\in\D$, the functions $k_{w_1},k_{w_2}\in\HH$ are linearly independent but not orthogonal. We will momentarily see that any $\HH$ relevant for Shimorin's question is irreducible. For a more general discussion of this technical condition, we refer to~\cite[Chapter~7]{AM02}.

Any reproducing kernel $K$ is \emph{positive}, in the sense that for any finite collection of points $w_1,\ldots,w_d\in\D$, the matrix $\begin{bmatrix}
	K(w_l,w_k)
\end{bmatrix}_{1\le k,l\le d}$ is positive semi-definite. We write $K\succeq0$ to indicate this property. It is an exercise to show that if a general function $K\colon\D\times\D\to\C$ is positive, then the functions $\{K(\cdot,w)\}_{w\in\D}$ generate a RKHS of functions on $\D$, see~\cite[Section~2.2]{AM02}. The space is uniquely determined by $K$, and consists of holomorphic functions if and only if $K$ is sesquiholomorphic. For any RKHS appearing in this paper, we assume that $K$ is sesquiholomorphic.

Assume $\HH$ is invariant under $M_z$ and $L$. If $T=(M_z,\HH)$, then $\ran T = \{f\in\HH\mid f(0)=0\}$. Indeed, if $f(0)=0$, then $f=TLf\in\ran T$. This establishes that $\{f\in\HH\mid f(0)=0\}\subset\ran T$. The reverse inclusion is trivial. From this we conclude that $\ran T = \{K^\HH(\cdot,0)\}^\perp$. In particular, $\ker T^*$ is one-dimensional.

A well-known example of a RKHS is the Hardy space $H^2$. This consists of functions $f\colon z\mapsto \sum_{n=0}^\infty f_nz^n$, where $\|f\|_{H^2}^2:=\sum_{n=0}^\infty|f_n|^2<\infty$. Its reproducing kernel
\[
K^{H^2}(z,w)=\frac{1}{1-z\overline{w}}=\sum_{n=0}^\infty \left(z\overline{w}\right)^n
\]
is the prototypical example of a normalized CNP kernel.

Consider an auxiliary Hilbert space $\HD$. The Schur class $\HS(\HD,\C)$ consists of holomorphic functions $B\colon \D\to \HD$ such that $\|B(z)\|_\HD < 1$ whenever $z\in\D$. Here we intentionally exclude constant functions with unit norm. It can be shown that if $B\in\HS(\HD,\C)$, then
\[
K^B(z,w)=\frac{1-\langle B(z),B(w)\rangle_\HD}{1-z\overline w} \succeq 0.
\]
The corresponding RKHS is called the \emph{de~Branges--Rovnyak space} $\HH[B]$. The correspondence between $\HH[B]$ and $K^B$ is one-to-one, although different Schur functions may define the same kernel. The \emph{rank} of $\HH[B]$ is the minimal dimension of the space $\HD$. Every $\HH[B]$-space is contractively contained in $H^2$, and the backward shift $L$ acts contractively on $\HH[B]$. If $\HH[B]$ is $M_z$-invariant, then $M_z$ is bounded and expansive on $\HH[B]$. We refer to \cite{AlemanMalman, BallBoloBasics, CGL24, debranges, FM162, LGR, sarasonbook} and the references therein for more background on de~Branges--Rovnyak spaces.

For any $\HH[B]$, if $w_1,w_2\in\D$ are distinct, then $K^B(w_1,w_2)\ne 0$, i.e., $k_{w_1}^B$ and $k_{w_2}^B$ are not orthogonal. If $f=\alpha_1 k_{w_1}^B + \alpha_2k_{w_2}^B$ and $\HH[B]$ is $M_z$-invariant, then $M_z^*f = \overline{w_1}\alpha_1 k_{w_1}^B + \overline{w_2}\alpha_2k_{w_2}^B$. Evaluating $f$ and $M_z^*f$ at any point shows that $\alpha_1,\alpha_2$ are uniquely determined by $f$. Hence, any $M_z$-invariant $\HH[B]$ is irreducible.

If $B\in\HS(\HD,\C)$ and $B(0)=0$, then we write $B\in\HS_0(\HD,\C)$. This corresponds precisely to the normalization $K^B(z,0)=1$.

If \(\mu\) is a finite positive Borel measure on \(\overline{\mathbb{D}}\), then
\[
U_\mu(w) = \int_{\partial\D} \frac{1 - |w|^2}{|w - z|^2} d\mu(z) + \int_{\mathbb{D}} \log \left| \frac{1 - \bar{w}z}{w - z} \right|^2 \frac{d\mu(z)}{1 - |z|^2}
\]
defines a positive superharmonic function on $\D$. Conversely, any non-trivial positive superharmonic function can be represented in this way. The Dirichlet type space \(D(\mu)\) consists of all holomorphic functions $f\colon\D\to\C$ such that
\begin{align*}
	\|f\|^2_{D(\mu)} := \|f\|_{H^2}^2 + \int_\D |f'(z)|^2 U_\mu(z)\,dA(z)<\infty.
\end{align*}
Here, $dA$ denotes area measure, normalized so that $A(\D)=1$. Dirichlet type spaces include the classical Dirichlet space. The above notion of $D(\mu)$ was introduced by Aleman~\cite{Al93}, following important works by Richter~\cite{Ri88,Ri91}. For more background, see also~\cite{ARSW19, EKMR14}. If $T=(M_z,D(\mu))$ and $\Delta = T^*T-I$, then
\begin{align}\label{eq:DeltaDmu}
	\langle\Delta f,f\rangle_{ D(\mu)}
	=\int_{\overline{\mathbb D}} |f(z)|^2\,d\mu(z).
\end{align}
The above identity appeared in \cite{Al93}. An easier to find reference is \cite[(5.2)]{LR15}. 

\section{Dirichlet type vs. concave shift}\label{sec:Comparison}

For $T=(M_z,D(\mu))$, \eqref{eq:DeltaDmu} easily implies that $T^*\Delta T\le \Delta$, so $T$ is concave. However, iteration of this argument shows that $T$ satisfies a much stronger condition. Define $\Delta_1=\Delta$, and $\Delta_{n}=T^*\Delta_{n-1} T-\Delta_{n-1}$ for $n\ge 2$. For each $n$, it then holds that
\begin{align*}
	\langle \Delta_{n} f,f\rangle_{D(\mu)}=(-1)^{n-1}\int_{\overline{\D}}|f(z)|^2(1-|z|^2)^{n-1}\, d\mu(z).
\end{align*}
Hence, the sequence $\Delta_1,\Delta_2,\Delta_3,\ldots$ has alternating signs. For $T$ to be concave, it would suffice that $\Delta_2\le 0$. 

To establish the existence of $\HH\in\Concave$ of non-Dirichlet type, we present a simple construction. Given $\zeta\in\overline{\D}$, consider the modified backward shift defined by $(L_\zeta f)(z)=\frac{f(z)-f(\zeta)}{z-\zeta}$. Let $w_1,\ldots,w_d$ be distinct points in $\D$, and define the joint evaluation map $E_{\vec{w}}\colon f\mapsto \begin{pmatrix} f(w_1) &\cdots & f(w_d) \end{pmatrix}^T$. For a positive definite matrix $C=\begin{bmatrix}
	c_{k,l}
\end{bmatrix}_{1\le k,l\le d}$, define the Hilbert space $\HH$ by imposing finiteness of the norm given by
\begin{align*}
	\|f\|_{\HH}^2=\|f\|_{H^2}^2 + \int_{\zeta\in\partial\D}\langle CE_{\vec{w}}L_\zeta f,E_{\vec{w}}L_\zeta f\rangle_{\C^d}\, dm(\zeta).
\end{align*}
Here, $dm$ denotes normalized arclength measure on $\partial\D$. $\HH$ is a RKHS, since it is contractively contained in $H^2$. It is clear that $\HH$ is in fact isomorphic to $H^2$, so $\HH$ is closed under $M_z$ and $L$. However, being isomorphic to $H^2$ is by no means sufficient for $\HH$ to be CNP, or even have concave forward shift.

It holds that $(L_\zeta M_z f)(z)=z(L_\zeta f)(z) + f(\zeta)$. By a calculation using this identity,
\begin{align*}
	\langle \Delta f,f\rangle_\HH = \langle C E_{\vec{w}}f,E_{\vec{w}}f\rangle_{\C^d}.
\end{align*}
Another straightforward calculation shows that $T$ is concave if and only if the matrix $\begin{bmatrix}
	c_{k,l}(1-\overline{w_k}w_l)
\end{bmatrix}_{1\le k,l\le d}$ is positive semi-definite, since this corresponds exactly to $\Delta \ge T^*\Delta T$. We also see that $\Delta f=0$ if and only if $E_{\vec{w}} f = 0$, so if $\HH$ is a $D(\mu)$-space, then \eqref{eq:DeltaDmu} implies that $\mu$ is supported on the points $w_1,\ldots,w_d$. The map $C\mapsto \HH$ is clearly injective, so if $C$ is not diagonal, then $\HH$ is not a $D(\mu)$-space.

\section{Logarithmically concave operators}\label{sec:ProofMainThm2}

In this section, we prove Theorem~\ref{mainthmonelogconnot} by constructing an example of a RKHS $\HH$, such that $T=(M_z,\HH)$ is logarithmically concave but $\HH$ is not CNP.

Let $S = (M_z, H^2)$ be the unilateral shift on $H^2$. Define the operator on $H^2\oplus\C^3$
	\[
	\mathcal T=
	\begin{pmatrix}S&\mathbf v\\0&A\end{pmatrix},
	\qquad \mathbf v=(4,-3,-1),\qquad
	A=\begin{pmatrix}
			0&5&-5/4\\
			0&0&6\\
			0&0&0
		\end{pmatrix}.
	\]
	Put
	\[
	p(z)=20z-3z^2,\qquad r(z)=4z^2,\qquad q(z)=120-23z-z^2,
	\qquad \tau=z\overline w.
	\]
	Set $G(z)=(r(z),p(z),q(z))$. Direct calculation gives
	\[
	zG(z)=z^3\mathbf v+G(z)A.
	\]
	Define
	\[
	U(f,x)(z)=\frac{z^3f(z)+G(z)x}{q(z)},
	\qquad f\in H^2,\quad x\in\C^3,
	\]
	and give its range the norm
	$\|U(f,x)\|_{\mathcal H}^2=\|f\|_{H^2}^2+\|x\|_{\C^3}^2$.
	The polynomials $r,p,q$ form a basis for the polynomials of degree at most two, so
	this decomposition is unique. Since $q$ and $1/q$ are bounded
	multipliers of $H^2$, the range is $H^2$ with an equivalent norm.
	So $\HH = H^2$, and $T$ and $L$ are bounded on $\HH$.
	
	The vectors
	\[
	(z^m,0)\quad(m\ge0),\qquad
	(0,e_1),(0,e_2),(0,e_3)
	\]
	form an orthonormal basis for $H^2 \oplus \C^3$.
	Their images under the map $U$ are
	\[
	\frac{z^{m+3}}{q(z)}\quad(m\ge0),\qquad
	\frac{r(z)}{q(z)},\quad
	\frac{p(z)}{q(z)},\quad 1.
	\]
	Since \(U\) is unitary, these functions form an orthonormal basis for $\HH$. Thus the reproducing kernel of $\HH$ is
	\[
	\begin{aligned}
			K(z,w)
			&=1+\frac{r(z)\overline{r(w)}+p(z)\overline{p(w)}}
			{q(z)\overline{q(w)}}+\sum_{m=0}^{\infty}
			\frac{z^{m+3}\overline w^{\,m+3}}
			{q(z)\overline{q(w)}}\\
			&=1+
			\frac{p(z)\overline{p(w)}+r(z)\overline{r(w)}
					+\dfrac{(z\overline w)^3}{1-z\overline w}}
			{q(z)\overline{q(w)}}\\
			& = 1+\frac{p(z)\overline{p(w)}+r(z)\overline{r(w)}
					+\dfrac{\tau^3}{1-\tau}}
			{q(z)\overline{q(w)}}.
		\end{aligned}
	\]
	Now we show $U\mathcal T=M_zU$.
	
	In fact, note that
	\begin{align*}
			&\mathcal T(f,x)=(zf+\mathbf vx,Ax),\\
			&zG(z)=z^3\mathbf v+G(z)A.
		\end{align*}
	We have
	\[
	\begin{aligned}
			(U\mathcal T(f,x))(z)
			&=\frac{z^3\bigl(zf(z)+\mathbf vx\bigr)+G(z)Ax}{q(z)}\\
			&=\frac{z^4f(z)+zG(z)x}{q(z)}\\
			&=z\,U(f,x)(z).
		\end{aligned}
	\]
	Thus
	$
	U\mathcal T=M_zU.
	$
	
	Next we show $T$ is expansive and logarithmically concave with $\rank \Delta = 3$. It is equivalent to establish that $\mathcal T$ has these properties.
	Note that $S^* \mathbf v = 0$. So
	\begin{align*}
			\HT^* \HT = \begin{pmatrix}
					I&0\\
					0&\mathbf v^* \mathbf v + A^*A\\
				\end{pmatrix}.
		\end{align*}
	Thus
	\begin{align*}
			\Delta = \HT^* \HT - I = \begin{pmatrix}
					0&0\\
					0&\mathbf v^* \mathbf v + A^*A - I\\
				\end{pmatrix},
		\end{align*}
	where
	\[
	D:= \mathbf v^* \mathbf v + A^*A - I
	=\begin{pmatrix}
			15&-12&-4\\
			-12&33&-13/4\\
			-4&-13/4&601/16
		\end{pmatrix}.
	\]
	Its leading principal minors are $15$, $351$, and $12186$.
	Therefore $\Delta \geq 0$ with $\rank \Delta = 3$.
	Note that
	\[
	\Delta+\Delta^2-\mathcal T^*\Delta\mathcal T
	=0\oplus E,
	\qquad E=D+D^2-A^*DA.
	\]
	Let
	\[
	B=\begin{pmatrix}
			4&-3&-1\\
			0&5&-5/4\\
			0&0&6
		\end{pmatrix}.
	\]
	Then $\mathbf v^* \mathbf v + A^*A = B^*B$ and $D = B^*B - I$.
	We have $A = PB$, where
	\[
	P=
	\begin{pmatrix}
			0&1&0\\
			0&0&1\\
			0&0&0
		\end{pmatrix}.
	\]
	Thus
	\[
	E=B^*WB,\qquad
	W= BB^*-I-P^*DP = \begin{pmatrix}
			25&-55/4&-6\\
			-55/4&169/16&9/2\\
			-6&9/2&2
		\end{pmatrix}.
	\]
	The leading principal minors of $W$ are
	\[
	25,\qquad 75,\qquad 6.
	\]
	Thus $W\geq0$ and $E\geq0$. This proves logarithmic concavity.
	
	Finally, we show that the reproducing kernel $K(z,w)$ is not a CNP kernel. Suppose to the contrary that \(K\) is a CNP kernel.
	Since $K(z,0)=1$, the complete Pick property would require
	$F=1-1/K$ to be a positive kernel.
	Then $Q(z,w)=q(z)\overline{q(w)}F(z,w)$ would also be positive.
	
	Write
	\[
	R(z,w)=p(z)\overline{p(w)}+r(z)\overline{r(w)}
	+\frac{\tau^3}{1-\tau}.
	\]
	Then
	\begin{align*}
			Q(z,w)&=q(z)\overline{q(w)}F(z,w) = \frac{q(z)\overline{q(w)} R(z,w)}{q(z)\overline{q(w)} + R(z,w)}\\
			& = \frac{R(z,w)}{1+[R(z,w)/q(z)\overline{q(w)}]}.
		\end{align*}
	Near zero,
	\[
	Q=R-\frac{R^2}{q(z)\overline{q(w)}}
	+\frac{R^3}{q(z)^2\overline{q(w)}^{\,2}}-\cdots.
	\]
	Only these three terms contribute to coefficients with both
	indices at most three. The coefficient matrix of $Q$ with
	indices $1,2,3$ is exactly
	\[
	C=\begin{pmatrix}
			400&-60&0\\[4pt]
			-60&125/9&65/54\\[4pt]
			0&65/54&125/432
		\end{pmatrix}.
	\]
	Note that every finite principal matrix of Taylor coefficients of a positive sesquiholomorphic kernel must be positive semi-definite.
	Indeed, suppose
	\[
	Q(z,w)=\sum_{i,j\ge0}c_{ij}z^i\overline w^{\,j}.
	\]
	Normalized derivative evaluation is a bounded functional on \(\mathcal H_Q\) (the Hilbert space determined by kernel $Q$), represented by
	\[
	g_j(z)=\frac1{j!}
	\left.\frac{\partial^j Q(z,w)}{\partial\overline w^{\,j}}\right|_{w=0}.
	\]
	Then
	\[
	c_{ij}
	=\frac1{i!\,j!}
	\left.
	\frac{\partial^{i+j}Q(z,w)}
	{\partial z^i\partial\overline w^{\,j}}
	\right|_{(0,0)}
	=\langle g_j,g_i\rangle_{\mathcal H_Q}.
	\]
	Thus for any finite index set \(I\) and scalars \((a_j)_{j\in I}\),
	\[
	\sum_{i,j\in I}\overline{a_i}c_{ij}a_j
	=\left\|\sum_{j\in I}a_jg_j\right\|_{\mathcal H_Q}^2
	\ge0.
	\]
	Since
	\[
	\det C=-\frac{10000}{729}<0,
	\]
	the matrix \(C\) is not positive semi-definite. Therefore \(Q\) cannot be a positive kernel. This contradiction proves that \(K\) is not a CNP kernel.

\section{Proof of Theorem~\ref{mainthmconcave}}\label{sec:ProofMainThm}

The goal in this section is to prove that if $\HH\in\Concave$, then $\HH$ is CNP. Our first lemma is a reduction to the case $\HH=\HH[B]$, where $B(0)=0$. The fundamental ingredient of the proof is that if $T=(M_z,\HH)$ is expansive, then $T$ is unitarily equivalent to some such $(M_z,\HH[B])$. Since unitary transformations in general don't preserve the CNP property, we also need to establish that the unitary transformation has a very particular form.
\begin{lemma}\label{hbmodsscnprhB}
	Let $\HH$ be a nonzero RKHS of holomorphic functions $f\colon\D\to\C$. Suppose $\HH$ is both $M_z$ and $L$ invariant. If $T = (M_z, \HH)$ is expansive, then there exists $B\in\HS_0\left(\overline{\ran\Delta},\C\right)$, and a nonvanishing holomorphic function $\delta\colon\D\to\C$ such that $K^{\HH}(z,w) = \delta(z) \overline{\delta(w)} K^B(z,w)$. In particular, $K^{\HH}$ is a CNP kernel if and only if $K^B$ is a CNP kernel. If $K^\HH$ satisfies the normalization condition $K^\HH(z,0)=1$, then $K^\HH=K^B$.
\end{lemma}
\begin{proof}
	It is clear that $T$ is analytic in the sense that $\bigcap_{n\ge 0}\ran T^n = \{0\}$. The fact that $T$ is unitarily equivalent to $(M_z,\HH[B])$ for some $B\in\HS_0\left(\overline{\ran\Delta},\C\right)$ follows from \cite[Theorem 4.6]{LGR}.
	
	Since $T$ is expansive, Fredholm theory implies that $\dim \ker (T^*-\lambda) = 1$ for $\lambda \in \D$. By \cite[Theorem 2.49]{AM02}, $K^\HH(z,w)=\delta(z) K^B(z, w)\overline{\delta(w)}$ for some nonvanishing holomorphic function $\delta\colon \D\to \C$. Since $K^{\HH}$ and $K^B$ define the same function $F$ via \eqref{eq:CNP}, $K^{\HH}$ is CNP if and only if $K^{B}$ is CNP.
	
	The condition $B(0)=0$ implies $K^{B}(z,0) = 1$. If $K^{\HH}(z,0) = 1$, then the identity $K^{\HH}(z,w) = \delta(z) \overline{\delta(w)} K^B(z,w)$ implies that $\delta$ is a unimodular constant.
\end{proof}

Recall that concave operators are expansive. In the light of Lemma~\ref{hbmodsscnprhB}, it's therefore sufficient to prove Theorem~\ref{mainthmconcave} under the stronger hypothesis $\HH\in\ConcaveNorm$. If $\HH\in\ConcaveNorm$ and $T=(M_z,\HH)$, then $\ker T^*$ is precisely the constant functions. It follows from \cite[Theorem~1]{Ri88} that polynomials are dense in $\HH$. Recall from Section~\ref{sec:Preliminaries} that any $M_z$-invariant $\HH[B]$-space is irreducible. The correspondence between $K^\HH$ and $K^B$ in Lemma~\ref{hbmodsscnprhB} implies that $\HH$ is irreducible.

We proceed with a simple structure theorem for the nilpotent case.
\begin{lemma}\label{agenmudrmf}
	If $\HH\in\ConcaveNil$, and $T=(M_z,\HH)$, then $T$ is unitarily equivalent to an operator of the form
	\begin{align}\label{modelopeatorconc}
		\begin{pmatrix}S&\beta\\0&A\end{pmatrix}
		\colon H^2\oplus\C^d\to H^2\oplus\C^d,
	\end{align}
	where $S = (M_z, H^2)$ is the Hardy shift, and $d=\rank \Delta$.
\end{lemma}
\begin{proof}
	This follows from the proof in \cite[Theorem 7.2]{LGR} or \cite[Lemma 3.2]{LR25}. Let $\HN = [\ran \Delta]_{T^*}$. By Remark~\ref{remark:NisranDelta}, $\HN=\ran\Delta$, so $\HM:=\HN^\perp$ is $M_z$-invariant. Moreover, $M_z|_\HM$ is isometric and analytic. Since $\dim \ran \Delta = d < \infty$, the subspace $\HM\ominus z\HM$ is 1-dimensional (\cite[Lemma 2.1]{ARS02}). Hence, $M_z|_\HM$ is unitarily equivalent to the unilateral shift $S = (M_z, H^2)$ of multiplicity 1.
\end{proof}

\begin{remark}
	The above lemma will be significantly improved in Theorem \ref{invariantsubspsigg}.
\end{remark}

Before proving the next lemma, we point out that if $A$ is nilpotent and $\psi$ is holomorphic in a neighbourhood of the origin, then $\psi(A)$ is well-defined using the standard holomorphic functional calculus.
\begin{lemma}\label{neltsdmpps}
	Let $\HH \in \mathfrak C_{\mathrm{nil}}$, $T=(M_z,\HH)$, $\Delta = T^*T-I$, $\HN=[\ran\Delta]_{T^*}$. If $\dim\HN=d$, then
	\[
	\mathcal{N} = \ker T^{*d}, \qquad \ker \Delta = z^d \mathcal{H}.
	\]
	If \(\psi\in\operatorname{Mult}(\mathcal{H})\), then
	\[
	\left. M_\psi^* \right|_{\mathcal{N}} = \psi(A)^*.
	\]
\end{lemma}
\begin{proof}
By Lemmas \ref{hbmodsscnprhB} and \ref{agenmudrmf}, $\HH = \HH[B]$ for some $B\in\HS_0(\ran \Delta,\C)$
and $T$ is unitarily equivalent to an operator of the form (\ref{modelopeatorconc}).
Note that $T^*|_\HN = A^*$ and $(A^*)^d = 0$. So $\HN \subseteq \ker T^{*d}$.
Since $\ker T^{*d} = (\ran T^d)^\perp$ is $d$-dimensional, $\HN = \ker T^{*d}$, and $\ker \Delta = \HN^\perp = z^d \mathcal{H}$.
	
	Suppose a multiplier $\psi$ has Taylor expansion
	\[
	\psi(z) = \sum_{n=0}^\infty a_n z^n = \sum_{n=0}^{d-1} a_n z^n + \psi_1(z).
	\]
	Let $f \in \HN, g \in \HH$. Then $\psi_1 g \in z^d \HH = \HN^\perp$. Thus
	\begin{align*}
		\langle M_\psi^* f, g\rangle& = \langle f, \psi g\rangle = \langle f, \sum_{n=0}^{d-1} a_n z^n g + \psi_1 g\rangle
		= \langle \sum_{n=0}^{d-1} \overline{a_n} T^{*n} f, g\rangle.
	\end{align*}
	Thus $M_\psi^*|_{\mathcal{N}} = \psi(A)^*$.
\end{proof}

The following result follows from von Neumann's inequality for operator-valued functions, see e.g. \cite[Theorem 3.42]{AMY20}.
We include a proof here for completeness.
\begin{lemma}\label{nilpotentestsumli}
	Let \(C\) be a nilpotent contraction on a finite-dimensional Hilbert space \(E\). If \(\psi = (\psi_\nu)_{\nu\geq 1}\) is an analytic row on \(\mathbb{D}\) with
	\[
	\sum_{\nu=1}^\infty |\psi_{\nu}(z)|^2 \leq 1 \qquad (z \in \mathbb{D}),
	\]
	then
	\[
	\sum_{\nu=1}^\infty \psi_\nu(C)\psi_\nu(C)^* \le I_E.
	\]
\end{lemma}
\begin{proof}
	Suppose $C^d=0$, and $\psi_\nu(z)=\sum_{k\ge0}a_{\nu,k}z^k$. Then $\psi_\nu(C) = \sum_{k=0}^{d-1} a_{\nu,k} C^k$. Fix $m\ge1$ and $0<r<1$, and let $\Psi_{m,r}(z) = (\psi_1(rz), \psi_2(rz), \ldots \psi_m(rz))$. $\Psi_{m,r}$ is analytic in a neighborhood of $\overline{\mathbb D}$ and satisfies
	\[
	\sup_{z\in\mathbb D}\|\Psi_{m,r}(z)\|
	=
	\sup_{z\in\mathbb D}
	\left(\sum_{\nu=1}^{m}|\psi_\nu(rz)|^2\right)^{1/2}
	\le1.
	\]
	By the matrix-valued von Neumann inequality, we have
	\[
	\left\|
	\begin{bmatrix}
		\psi_1(rC)&\cdots&\psi_m(rC)
	\end{bmatrix}
	\right\|
	\le1.
	\]
	Consequently,
	\[
	\sum_{\nu=1}^{m}
	\psi_\nu(rC)\psi_\nu(rC)^*\le I_E.
	\]
	Since
	\[
	\psi_\nu(rC)
	=
	\sum_{k=0}^{d-1}a_{\nu,k}r^kC^k
	\longrightarrow
	\sum_{k=0}^{d-1}a_{\nu,k}C^k
	=
	\psi_\nu(C)
	\]
	in operator norm as $r\uparrow1$, for every $m\ge1$,
	\[
	\sum_{\nu=1}^{m}
	\psi_\nu(C)\psi_\nu(C)^*\le I_E.
	\]
	It follows that
	\[
	\sum_{\nu=1}^{\infty}
	\psi_\nu(C)\psi_\nu(C)^*\le I_E.
	\]
\end{proof}

We need the following rank-one extension scheme introduced by Stefan Richter. It was also used in \cite[Section 2.2]{GL}.
\begin{lemma}\label{rankoneextens}
	Let $\HH \in \ConcaveNil$. Assume that $\HH$ is CNP, and that
	\begin{align*}
		\widetilde{T} = \begin{pmatrix}T&u \otimes 1\\0&0\end{pmatrix}
		\colon \HH\oplus\C\to\HH\oplus\C
	\end{align*}
	is concave for some $u\in\HH$ with $u(0) \neq 0$. Then $u$ has no zero on $\D$. If $g = z/u$, and
	\[
	\widetilde{\mathcal{H}} := \{a + gf : f \in \mathcal{H}, \ a \in \mathbb{C}\}, \qquad
	\|a + gf\|_{\widetilde{\mathcal{H}}}^2 := |a|^2 + \|f\|_{\mathcal{H}}^2,
	\]
	then $\widetilde{\HH} \in \ConcaveNorm$ and $(M_z, \widetilde{\HH})$ is unitarily equivalent to $\widetilde{T}$. Moreover, $\widetilde{\HH}$ is CNP, with kernel
	\begin{align}\label{khtrelkhg}
		K^{\widetilde{\HH}}(z,w) = 1 + g(z) \overline{g(w)} K^\HH(z,w).
	\end{align}
\end{lemma}
\begin{proof}
	Let \(P = T^*T = I + \Delta\) and \(P_0 = 1 \otimes 1\) on \(\mathcal{H}\). We have $1 \perp z \HH$ and \(\|1\| = 1\).
	Note that
	\begin{equation*}
		L = (T^*T)^{-1} T^* = P^{-1}T^*, \qquad LL^* = P^{-1}, \qquad TP^{-1}T^* = I - P_0.
	\end{equation*}
	
	Since $K^\HH$ is a CNP kernel, there are analytic functions $c_\nu$ on $\D$ such that
	\[
	F(z, w) = 1 - \frac{1}{K^\HH(z,w)} = \sum_{\nu} c_\nu(z) \overline{c_\nu(w)}, \qquad c_\nu(0) = 0,
	\]
	see \cite{AM02}. The identity \((1 - F)K^\HH = 1\) shows that
	\begin{equation*}
		\sum_{\nu} M_{c_\nu} M_{c_\nu}^* = I - P_0.
	\end{equation*}
	Let $\psi_\nu = Lc_\nu$, so that $M_{\psi_\nu} = L M_{c_\nu}$, and define
	\[
	\Phi(X) = \sum_{\nu} M_{\psi_\nu} X M_{\psi_\nu}^*.
	\]
	Then $\Phi$ is a bounded completely positive map. Since $L1 = 0$, we have
	\begin{equation*}
		\Phi(I) = L \sum_{\nu} M_{c_\nu} M_{c_\nu}^* L^* = L(I - P_0)L^* = P^{-1}.
	\end{equation*}
	
	By Lemma~\ref{hbmodsscnprhB}, $\HH = \HH[B]$ for $B \in \HS_0(\ran \Delta,\C)$. Since $\HH[B]$ is contractively contained in $H^2$ (\cite{BallBoloBasics}),
	we have $K^\HH(z,z) \leq (1-|z|^2)^{-1}$ (\cite{PaulsenRag}). It follows that if $z\ne0$, then
	\[
	\sum_{\nu} |\psi_{\nu}(z)|^{2} = \frac{1 - 1/K^\HH(z,z)}{|z|^{2}} \leq 1,
	\]
	and this inequality extends to $z=0$ by continuity.
	
	We claim that
	\begin{align}\label{phidelest}
		\Phi(\Delta) \geq P^{-1}\Delta P^{-1}.
	\end{align}
	If $\Delta = 0$, then it holds. Otherwise, let $D = \Delta|_{\mathcal{N}}$, $C= D^{1/2}AD^{-1/2}$, where $A=P_{\HN}T|_{\HN}$, and let $\iota : \mathcal{N} \to \mathcal{H}$ be inclusion. Then $C$ is nilpotent, and since $T$ is concave, $C\colon \HN\to \HN$ is contractive. Lemma \ref{nilpotentestsumli} implies that $\sum_{\nu} \psi_{\nu}(C)\psi_{\nu}(C)^{*} \leq I$. We rephrase this by saying that the column operator
	\[
	\Gamma = \big(D^{-1/2}\psi_{\nu}(A)^{*}D^{1/2}\big)_{\nu} = \big(\psi_{\nu}(C)^{*}\big)_{\nu}
	\]
	is contractive. In tandem with $\Gamma$, we consider the row operator $W=\big( M_{\psi_{\nu}}\iota D^{1/2}\big)_\nu$. By calculation,
	\begin{align*}
		WW^*
		=
		\sum_{\nu} M_{\psi_\nu} \iota D\iota^* M_{\psi_\nu}^*
		=
		\sum_{\nu} M_{\psi_\nu} \Delta M_{\psi_\nu}^*
		=
		\Phi(\Delta)
		,
	\end{align*}
	and
	\begin{align*}
		W\Gamma = \sum_{\nu} M_{\psi_{\nu}} \iota \psi_{\nu}(A)^* D^{1/2} = \sum_{\nu} M_{\psi_{\nu}} M_{\psi_{\nu}}^* \iota D^{1/2} = P^{-1} \iota D^{1/2}.
	\end{align*}
	In the latter calculation we have used Lemma~\ref{neltsdmpps}. Since $\Gamma$ is contractive, $\Phi(\Delta)=WW^*\ge W\Gamma\Gamma^*W^* = P^{-1}\Delta P^{-1}$, thus proving (\ref{phidelest}).
	
	We now isolate the first row of $\widetilde{T}$, by letting
	\[
	X = (T \ \ u \otimes 1) : \HH \oplus \mathbb{C} \to \HH.
	\]
	Clearly $X$ is bounded. Moreover,
	\[
	X\begin{pmatrix}
		f\\ a
	\end{pmatrix} = h \iff
	\begin{pmatrix}
		f \\ a
	\end{pmatrix}
	=\begin{pmatrix}
		L\big(h-uh(0)/u(0)\big) \\ h(0)/u(0)
	\end{pmatrix},
	\]
	so $X$ is invertible.
	
	Let \(\widetilde{\Delta} = \widetilde{T}^*\widetilde{T} - I\). One checks that
	\begin{align*}
		\widetilde{\Delta} = X^* X - I, \qquad \widetilde{T}^* \widetilde{\Delta} \widetilde{T} = X^* \Delta X.
	\end{align*}
	Hence,
	\begin{equation*}
		\widetilde{\Delta} - \widetilde{T}^*\widetilde{\Delta}\widetilde{T} = X^* (I - \Delta) X - I \geq 0.
	\end{equation*}
	In particular, \(Y = I - \Delta\) is strictly positive and boundedly invertible. Since $ Y \geq (X^*)^{-1} X^{-1}$,
	it follows that
	\begin{equation}\label{Tsttputu}
		TT^* + u \otimes u = XX^* \geq Y^{-1}.
	\end{equation}
	On the other hand, $\Phi(I) = P^{-1}$, and (\ref{phidelest}) gives
	\begin{equation*}
		\Phi(Y) = P^{-1} - \Phi(\Delta) \leq P^{-1} - P^{-1}\Delta P^{-1} = P^{-1} - P^{-1}(P-I) P^{-1} = P^{-2}.
	\end{equation*}
	Since $\Phi$ is completely positive, applying \(\Phi\) entrywise to the positive operator matrix
	\[
	\begin{pmatrix}
		Y^{-1} & I \\
		I & Y
	\end{pmatrix}
	\geq 0,
	\]
	we obtain that
	\[
	\begin{pmatrix}
		\Phi(Y^{-1}) & P^{-1} \\
		P^{-1} & P^{-2}
	\end{pmatrix}
	\geq \begin{pmatrix}
		\Phi(Y^{-1}) & \Phi(I) \\
		\Phi(I) & \Phi(Y)
	\end{pmatrix}\geq
	0.
	\]
	The left-hand side has Schur complement $\Phi(Y^{-1}) - P^{-1}(P^{-2})^{-1}P^{-1}$, and this is positive, so
	\[
	\Phi(Y^{-1}) \geq P^{-1}(P^{-2})^{-1}P^{-1} = I.
	\]
	Since all $M_{\psi_{\nu}}$ commute with $T$,
	\[
	\Phi(TT^*) = T\Phi(I)T^* = TP^{-1}T^* = I - P_0.
	\]
	Applying $\Phi$ to (\ref{Tsttputu}) and using $\Phi(Y^{-1}) \geq I$ gives
	\[
	\Phi(u \otimes u) \geq P_0.
	\]
	
	Note that
	\[
	\Phi(u \otimes u) = \sum_\nu \psi_\nu u \otimes \psi_\nu u.
	\]
	Thus
	\begin{align}\label{poskelufmo}
		0& \preceq\langle(\Phi(u \otimes u) - P_0) k^\HH_w, k^\HH_z\rangle \notag\\
		& = \sum_\nu \psi_\nu(z) u(z)\overline{\psi_\nu(w) u(w)} - 1 \notag\\
		& = u(z)\overline{u(w)}\frac{F(z,w)}{z\overline{w}} - 1.
	\end{align}
	In particular, $|u(z)|^2 \frac{F(z,z)}{|z|^2} - 1 \geq 0$. We have already seen that $F(z,z)\le |z|^2$. Therefore, $u$ has no zeros in $\mathbb{D}$.
	
	Recall that $g = z/u$. Then (\ref{poskelufmo}) implies that
	\[
	G(z,w) := F(z,w) - g(z)\overline{g(w)} \succeq 0.
	\]
	Moreover,
	\[
	0 \leq G(z,z) \leq F(z,z) = 1 - \frac{1}{K^\HH(z,z)} < 1.
	\]
	Hence $|G(z,w)| < 1$. The kernel in (\ref{khtrelkhg}) satisfies $K^{\widetilde{\HH}} = K^\HH (1 - G)$ and
	\[
	1 - \frac{1}{K^{\widetilde{\HH}}(z,w)} = \frac{g(z) \overline{g(w)}K^\HH(z,w)}{K^{\widetilde{\HH}}(z,w)} = \frac{g(z) \overline{g(w)}}{1 - G(z,w)} = g(z) \overline{g(w)} \sum_{m=0}^{\infty} G^m \succeq 0.
	\]
	This shows that $K^{\widetilde{\HH}}$ is a CNP kernel.
	
	Now we show $\widetilde{\HH} \in \ConcaveNorm$ and $(M_z, \widetilde{\HH})$ is unitarily equivalent to $\widetilde{T}$.
	Since $\C \cap g \HH = \{0\}$, we have $\widetilde{\HH} = g \HH\oplus \C $ and $1 + g(z) \overline{g(w)} K^\HH(z,w)$ is the kernel for $\widetilde{\HH}$ (\cite{PaulsenRag}).
	Let
	\[
	U : \HH \oplus \mathbb{C} \to \widetilde{ \HH}, \qquad U(f, a) = a + \frac{zf}{u} = a + gf.
	\]
	Then $U$ is unitary and
	\[
	U\widetilde{T}(f,a) = \frac{z(Tf + a u)}{u} = z U(f,a).
	\]
	Thus $(M_z, \widetilde{\HH})$ is unitarily equivalent to $\widetilde{T}$. Let $h = a + gf \in \widetilde{\HH}$. Then
	\begin{align*}
		zh = az + zgf = g(au+Tf) = g X(f,a).
	\end{align*}
	Since $X = (T \ \ u \otimes 1) : \HH \oplus \mathbb{C} \longrightarrow \HH$ is onto,
	we have $z \widetilde{\HH} = g\HH$ and $\widetilde{\HH} \ominus z \widetilde{\HH} = \C$. So polynomials are dense in $\widetilde{\HH}$ (\cite[Theorem 1]{Ri88})
	and $\widetilde{\HH} \in \ConcaveNorm$.
\end{proof}

We are now ready to prove the following preliminary version of Theorem~\ref{mainthmconcave}.
\begin{proposition}\label{nilpotconcnp}
	Let $\HH \in \mathfrak C_{\mathrm{nil}}$. Then $\HH$ is a CNP space.
\end{proposition}
\begin{proof}
	By Lemma \ref{neltsdmpps}, $\ker \Delta = z^d \mathcal{H}$, where $d = \dim \ran \Delta$. Let
	\begin{align*}
		\mathcal{M}_j = z^{d-j}\mathcal{H} \quad (0 \le j \le d), \quad T_j = T|_{\mathcal{M}_j}.
	\end{align*}
	For $j \ge 1$ we have
	\begin{align*}
		z\mathcal{M}_j = \mathcal{M}_{j-1}, \quad \dim(\mathcal{M}_j \ominus \mathcal{M}_{j-1}) = 1.
	\end{align*}
	Every $T_j$ is expansive and concave, with defect
	\begin{align*}
		\Delta_j = T_j^*T_j-I=P_{\mathcal{M}_j}\Delta|_{\mathcal{M}_j}.
	\end{align*}
	Positivity of $\Delta$ gives
	\begin{align}\label{kendeljnj}
		\ker \Delta_j = \mathcal{M}_0, \quad \mathcal{N}_j := \operatorname{ran}\Delta_j = \mathcal{M}_j \ominus \mathcal{M}_0, \quad \operatorname{rank}\Delta_j = j.
	\end{align}
	Furthermore, $T_j^j \mathcal{M}_j = \mathcal{M}_0$, so its finite compression satisfies
	\begin{align}\label{nilpnjc}
		\big(P_{\mathcal{N}_j}T_j|_{\mathcal{N}_j}\big)^j = 0 \quad (j \ge 1).
	\end{align}
	
	If $d = 0$, then $\HH = H^2$, and our conclusion holds.
	
	For each \(1 \leq j \leq d\), choose unit vectors $\eta_j \in \mathcal{M}_j \ominus \mathcal{M}_{j-1}$. We inductively construct CNP spaces \(\mathcal{K}_j \in \mathfrak C_{\mathrm{nil}}\), and unitaries
	\begin{equation}\label{indwjgbmej}
		W_j : \mathcal{K}_j \longrightarrow \mathcal{M}_j, \qquad W_j h = \eta_j h,
	\end{equation}
	intertwining multiplication by \(z\). We let $\eta_d = 1$.
	
	The initial construction supplies \(j = 0\). In fact, by Lemma \ref{agenmudrmf}, $T|_{\HM_0}$ is unitarily equivalent to $(M_z, H^2)$.
	So take $\HK_0 = H^2$, $W_0 f = \eta_0 f$, where $\eta_0 \in \HM_0 \ominus z \HM_0$ is a unit vector.
	
	Suppose \(j < d\) and that we have constructed $\HK_j$ and $W_j$. Because \(T\eta_{j+1} \in \mathcal{M}_j\), set
	\[
	u_j = W_j^* T\eta_{j+1} \in \mathcal{K}_j.
	\]
	We have \(u_j(0) \neq 0\). Otherwise, boundedness of $L$ on \(\mathcal{K}_j\) would give \(u_j \in z\mathcal{K}_j\), and hence \(T\eta_{j+1} \in T\mathcal{M}_j\). Injectivity of \(T\) would imply \(\eta_{j+1} \in \mathcal{M}_j\), a contradiction. Let $V_j: \mathcal{K}_j \oplus \mathbb{C} \rightarrow \HM_{j+1} = \mathcal{M}_j \oplus \mathbb{C}\eta_{j+1}$ be defined by
	\[
	V_j(f,a) = W_j f + a \eta_{j+1}.
	\]
	Then $V_j$ is unitary and
	\begin{align*}
		T_{j+1} V_j (f,a) = z (W_j f + a \eta_{j+1}) = W_j(zf + au_j)=V_j(zf + au_j, 0).
	\end{align*}
	Thus \(T_{j+1}\) is unitarily equivalent to
	\[
	\begin{pmatrix}
		M_z & u_j \\
		0 & 0
	\end{pmatrix}
	\colon \mathcal{K}_j \oplus \mathbb{C}\to \mathcal{K}_j \oplus \mathbb{C}.
	\]
	It is concave since \(T_{j+1}\) is concave. By Lemma \ref{rankoneextens}, $u_j$ has no zero on $\D$. Let
	\[
	\mathcal{K}_{j+1} = \left\{ a + \frac{z f}{u_j} : f \in \mathcal{K}_j,\ a \in \mathbb{C} \right\}
	\]
	with squared norm $|a|^2 + \|f\|_{\mathcal{K}_j}^2$. By Lemma \ref{rankoneextens}, \(\mathcal{K}_{j+1} \in \mathfrak C_{\mathrm{nil}}\) is a CNP space.
	Define
	\[
	W_{j+1}\left(a + \frac{z f}{u_j}\right) = a \eta_{j+1} + W_j f.
	\]
	This is a unitary onto $\mathcal{M}_{j+1}$. The identity $z\eta_{j+1} = W_j u_j = \eta_j u_j$ shows pointwise that
	\[
	a\eta_{j+1} + \eta_j f = \eta_{j+1}\left(a + \frac{z f}{u_j}\right).
	\]
	Thus (\ref{indwjgbmej}) holds at $j+1$, and this unitary intertwines multiplication by $z$. The induction is complete.
	
	Finally, $\mathcal{M}_d = \mathcal{H}$ and $\eta_d = 1$. Consequently $W_d h = h$ pointwise, so $\mathcal{K}_d$ and $\mathcal{H}$ are the same Hilbert function space with the same kernel. Since $\mathcal{K}_d$ is CNP, so is $\mathcal{H}$.
\end{proof}

For the next lemma, recall Definition~\ref{def:G}.

\begin{lemma}\label{graphnormdef}
	Let $E$ be a Hilbert space, $C\in\HB(E)$ a contraction, and $v\in E$. The space \(\mathcal{G}=\mathcal{G}(E,C,v)\) is a RKHS of holomorphic functions on $\D$. Polynomials form a dense subspace of $\mathcal{G}$. If \(p\) is a polynomial, then
	\begin{equation}\label{solutionfppc}
		y(p)_k = (L^k p)(C)v = \sum_{i \geq 0} p_{k+i} C^i v,
	\end{equation}
	where the sum is finite. Moreover, \(K^{\HG}(z, 0) = 1\).
\end{lemma}
\begin{proof}
	We claim that the set of pairs $(f,y) \in H^2 \oplus \ell^2(\mathbb{N},E)$ satisfying (\ref{fysaiedn}) is a closed linear subspace. Indeed, if $f^{(j)}\to f$ in $H^2$, and $y^{(j)}=y(f^{(j)})\to y$ in $\ell^2(\mathbb{N},E)$, then
	\[
	f_k^{(j)}\to f_k,\qquad
	y^{(j)}_{k+1}\to y_{k+1},
	\]
	for each $k\ge 0$. Passing to the limit in \eqref{fysaiedn} yields $y_k-Cy_{k+1}=f_kv$, so $y=y(f)$. This closed subspace is isometrically isomorphic with $\mathcal{G}$, so $\mathcal{G}$ is complete. Being contractively contained in $H^2$, $\mathcal{G}$ is a RKHS. By calculation, formula (\ref{solutionfppc}) gives a square-summable solution for every polynomial.
	
	We now prove polynomials are dense in $\mathcal{G}$. Let $f \in \mathcal{G}$, put $p_m(z) = \sum_{j=0}^m f_j z^j$.
	Iterating (\ref{fysaiedn}) gives
	\[
	y_k = \sum_{i=0}^{m-k} f_{k+i} C^i v + C^{m-k+1} y_{m+1} \qquad (1 \le k \le m).
	\]
	Consequently
	\[
	y(p_m)_k -y_k=
	\begin{cases}
		- C^{m-k+1} y_{m+1}, & 1 \le k \le m,\\
		-y_k, & k > m.
	\end{cases}
	\]
	It follows that
	\[
	\|y(p_m) - y\|_{\ell^2(E)}^2 \le m\|y_{m+1}\|_E^2 + \sum_{k>m} \|y_k\|_E^2.
	\]
	Now note that $\liminf_{m\to\infty}m\|y_{m+1}\|^2=0$, for otherwise there would exist $\varepsilon>0$ and $M$ such
	that $\|y_{m+1}\|^2\ge\varepsilon/m$ for every $m\ge M$. This would immediately contradict square summability. It follows that $m_j\|y_{m_j+1}\|^2\to0$ for some subsequence. Then
	\[
	\begin{aligned}
		\|p_{m_j}-f\|_{\mathcal{G}}^2
		&=\|p_{m_j}-f\|_{H^2}^2
		+\|y(p_{m_j}) - y\|_{\ell^2(E)}^2\\
		&\le
		\sum_{n>m_j}|f_n|^2
		+m_j\|y_{m_j+1}\|_E^2
		+\sum_{k>m_j}\|y_k\|_E^2
		\longrightarrow0.
	\end{aligned}
	\]
	This proves that polynomials are dense in $\mathcal{G}$.
	
	Finally, \(y(1) = 0\), so \(\|1\|_\HG  = 1\) and \(\langle f, 1 \rangle_{\mathcal{G}} = f(0)\). This proves \(K^\HG(z, 0) = 1\).
\end{proof}

\begin{proposition}\label{frknilpappx}
	Let $\HG=\mathcal{G}(E,C,v)$ be as in Lemma \ref{graphnormdef}. There exist $\HH_1, \HH_2,\ldots \in \mathfrak C_{\mathrm{nil}}$ with the following properties. If $T_n := (M_z, \HH_n)$, $\Delta_n := T_n^*T_n - I$, $\HN_n := \ran\Delta_n$, $A_n := P_{\HN_n} T_n|_{\HN_n}$, then $\rank \Delta_n \leq n$, and $K^{\HH_n} \to K^\HG$ pointwise on $\D\times\D$. Consequently, $\HG$ is CNP.
\end{proposition}
\begin{proof}
	Let $\{e_i: 0 \leq i \leq n-1\}$ be an orthonormal basis for $\C^n$ and let \(S_n\) be the truncated forward shift on \(\mathbb{C}^n\):
	\[
	S_n e_j = e_{j+1} \quad (0 \leq j < n - 1), \qquad S_n e_{n-1} = 0.
	\]
	Set
	\[
	h_n = \frac{1}{\sqrt{n}} \sum_{j=0}^{n-1} e_j, \qquad C_n = C \otimes S_n, \qquad v_n = v \otimes h_n.
	\]
	Since $\|S_n\| \leq 1$ and $S_n^n = 0$, we have \(\|C_n\| \leq 1\) and \(C_n^n = 0\).
	
	For \(f \in H^2\) define
	\begin{equation}\label{defnmihnc}
		\|f\|_n^2 = \|f\|_{H^2}^2 + \sum_{k \geq 1} \|(L^k f)(C_n)v_n\|_{E\otimes \C^n}^2.
	\end{equation}
	By the nilpotency of $C_n$,
	\[
	(L^k f)(C_n)v_n = \sum_{i=0}^{n-1} f_{k+i} C_n^i v_n.
	\]
	For each \(j \geq 0\), consider the operator \(Q_j : H^2 \to \ell^2(\mathbb{N}, E)\) given by
	\[
	(Q_j f)_k = \sum_{i=0}^{j} f_{k+i} C^i v \qquad (k \geq 1).
	\]
	The triangle inequality gives \(\|Q_j\| \leq (j + 1)\|v\|_E\). Note that
	\begin{align*}
		C_n^i v_n &= C^i v \otimes S_n^i h_n = C^i v \otimes \frac{1}{\sqrt{n}} \sum_{j=0}^{n-1-i} e_{j+i}
		= \frac{1}{\sqrt{n}} \sum_{j=i}^{n-1} C^i v \otimes e_{j}.
	\end{align*}
	Thus
	\begin{align*}
		(L^k f)(C_n)v_n & = \sum_{i=0}^{n-1} f_{k+i} \frac{1}{\sqrt{n}} \sum_{j=i}^{n-1} C^i v \otimes e_{j}\\
		& = \frac{1}{\sqrt{n}} \sum_{j=0}^{n-1} \sum_{i=0}^{j} f_{k+i} C^i v \otimes e_j\\
		& =  \frac{1}{\sqrt{n}} \sum_{j=0}^{n-1} (Q_j f)_k \otimes e_j.
	\end{align*}
	So
	\begin{align}\label{anoexpfnom}
		\|f\|_n^2 = \|f\|_{H^2}^2 + \frac1n \sum_{j=0}^{n-1} \|Q_j f\|_{\ell^2(E)}^2.
	\end{align}
	In particular,
	\[
	\|f\|_{H^2}^2 \le \|f\|_n^2 \le \big(1 + n^2\|v\|^2\big)\|f\|_{H^2}^2.
	\]
	Thus this is an equivalent Hardy norm, and we denote this space by $\mathcal{H}_n$.
	$(M_z, \HH_n)$ and $L$ are bounded and polynomials are dense.
	Since every $Q_j 1 = 0$, the constant $1$ is a unit vector representing evaluation at zero, so $K^{\HH_n}(z,0)=1$.
	
	Define
	\[
	E_n = \operatorname{span}\{C_n^j v_n : 0 \le j < n\}, \quad
	J_n f = f(C_n) v_n = \sum_{i=0}^{n-1} f_i C_n^i v_n, ~ f \in \HH_n.
	\]
	The map $J_n : \mathcal{H}_n \to E_n$ is bounded and onto. Notice that $\dim E_n \le n$, independently of the dimension of $E$.
	By (\ref{defnmihnc}), we have
	\[
	\|z f\|_n^2 - \|f\|_n^2 = \|J_n f\|^2,\quad J_n T_n = C_n J_n.
	\]
	It follows that
	\[
	\Delta_n = J_n^* J_n,\quad \Delta_n - T_n^* \Delta_n T_n = J_n^*(I - C_n^* C_n) J_n \ge 0.
	\]
	Thus \(T_n\) is expansive and concave, with \(\operatorname{rank} \Delta_n = \dim E_n \leq n\).
	
	Recall that $\HN_n = \ran \Delta_n$. Then $\HN_n = \ran J_n^*$ and \(\mathcal{N}_n^\perp = \ker J_n\). So \(J_n P_{\mathcal{N}_n} = J_n (P_{\mathcal{N}_n} + P_{\mathcal{N}_n^\perp}) = J_n\). The restriction \(J_n|_{\mathcal{N}_n} : \mathcal{N}_n \to E_n\) is an isomorphism and
	\[
	(J_n|_{\mathcal{N}_n})A_n = (C_n|_{E_n})(J_n|_{\mathcal{N}_n}).
	\]
	Consequently \(A_n\) is similar to \(C_n|_{E_n}\), proving \(A_n^n = 0\). If \(E_n = \{0\}\), the same assertion is understood on the zero space.
	This proves that $\mathcal{H}_n \in \mathfrak C_{\mathrm{nil}}$.
	
	We are left to show the kernel convergence. Fix \(w \in \mathbb{D}\) and set
	\[
	u_n = K^{\HH_n}(\cdot, w), \qquad \bar{y}_n = \frac{1}{n} \sum_{j=0}^{n-1} Q_j u_n.
	\]
	The Cauchy--Schwarz inequality and (\ref{anoexpfnom}) imply
	\begin{equation}\label{normbduhke}
		\|u_n\|_{H^2}^2 + \|\bar{y}_n\|_{\ell^2(E)}^2 \leq \|u_n\|_n^2 = K^{\HH_n}(w, w) \leq \frac{1}{1 - |w|^2}.
	\end{equation}
	The last inequality follows from \(\|f\|_{H^2} \leq \|f\|_n\) and the Hardy evaluation bound.
	Hence, every subsequence has a further weakly convergent subsequence (which we continue to denote using the same indices),
	\[
	u_n \rightarrow u \text{ weakly in } H^2, \qquad \bar{y}_n \rightarrow y \text{ weakly in } \ell^2(\mathbb{N}, E).
	\]
	
	For each $k \ge 1$,
	\begin{align*}
		(Q_ju_n)_k = \sum_{i=0}^{j} (u_n)_{k+i} C^i v, \qquad C(Q_ju_n)_{k+1} = \sum_{i=0}^{j} (u_n)_{k+1+i} C^{i+1} v.
	\end{align*}
	So
	\begin{align*}
		(Q_ju_n)_k - C(Q_ju_n)_{k+1} = (u_n)_k v - (u_n)_{k+j+1} C^{j+1} v.
	\end{align*}
	Then
	\[
	(\bar y_n)_k - C(\bar y_n)_{k+1} = (u_n)_k v - \varepsilon_{n,k}, \quad
	\varepsilon_{n,k} = \frac1n \sum_{j=0}^{n-1} (u_n)_{k+j+1} C^{j+1} v.
	\]
	By Cauchy--Schwarz,
	\[
	\|\varepsilon_{n,k}\| \le \frac{\|v\|}{\sqrt{n}} \|u_n\|_{H^2} \longrightarrow 0.
	\]
	Coefficient evaluation on $H^2$, coordinate evaluation on $\ell^2(E)$, and $C$ are bounded maps. Taking weak limits therefore gives $y_k - C y_{k+1} = u_k v$ for every $k$. Lemma \ref{graphnormdef} then implies $u \in \mathcal{G}$ and $y = y(u)$.
	
	Let $p$ be a polynomial of degree $d$. By (\ref{solutionfppc}), $Q_j p = y(p)$ whenever $j \ge d$.
	The reproducing identity and polarization of (\ref{anoexpfnom}) give
	\[
	\begin{aligned}
		p(w) &= \langle p, u_n\rangle_{H^2} + \frac1n \sum_{j=0}^{n-1} \langle Q_j p, Q_j u_n\rangle_{\ell^2(E)} \\
		&= \langle p, u_n\rangle_{H^2} + \langle y(p), \bar y_n\rangle_{\ell^2(E)} + \frac1n \sum_{j=0}^{d-1} \langle Q_j p - y(p), Q_j u_n\rangle_{\ell^2(E)}.
	\end{aligned}
	\]
	For $n > d$, the absolute value of $\frac1n \sum_{j=0}^{n-1} \langle Q_j p - y(p), Q_j u_n\rangle_{\ell^2(E)}$ is bounded by
	\[
	\left(\frac1n \sum_{j=0}^{d-1} \|Q_j p - y(p)\|_{\ell^2(E)}^2\right)^{1/2}
	\left(\frac1n \sum_{j=0}^{n-1} \|Q_j u_n\|_{\ell^2(E)}^2\right)^{1/2}.
	\]
	The first factor is $O(n^{-1/2})$, and the second is bounded by (\ref{normbduhke}). Taking limits then gives
	\[
	p(w) = \langle p, u\rangle_{H^2} + \langle y(p), y(u)\rangle_{\ell^2(E)} = \langle p, u\rangle_{\mathcal{G}}.
	\]
	Since polynomials are dense in $\mathcal{G}$, we conclude that $u = K^\HG(\cdot,w)$.
	Every subsequence has this same weak cluster point, so the full sequence $u_n = K^{\HH_n}(\cdot, w)$ converges weakly to $K^\HG (\cdot,w)$ in $H^2$.
	Applying bounded point evaluation on $H^2$, we obtain that $K^{\HH_n}(z,w) \to K^\HG (z,w)$ for $z,w\in\mathbb{D}$. It follows from \cite[Lemma 2.2]{Sh02} that $\HG$ is CNP.
\end{proof}

\begin{proof}[Proof of Proposition~\ref{prop:Gspaces} and Theorem~\ref{mainthmconcave}]
	Let
	\[
	E = \overline{\ran \Delta^{1/2}} = \overline{\ran \Delta}, \qquad v = \Delta^{1/2}1.
	\]
	Define $C: \ran \Delta^{1/2} \rightarrow E$ by
	\begin{equation}\label{defnoccvigs}
		C(\Delta^{1/2}f) = \Delta^{1/2}Tf.
	\end{equation}
	Since $T$ is concave,
	\[
	\|\Delta^{1/2}Tf\|^2 \leq \|\Delta^{1/2}f\|^2,
	\]
	which implies that $C$ is well defined and contractive on \(\ran \Delta^{1/2}\). It extends to a contraction \(C\) on \(E\).
	For every polynomial \(p\), (\ref{defnoccvigs}) gives
	\[
	\Delta^{1/2}p = p(C)v.
	\]
	Since $\HH = \HH[B]$ for some $B \in \HS_0(E,\C)$, we have $1 \perp z\HH$ and $\|1\|_\HH = 1$.
	By \(p = p(0) + TLp\),
	\[
	\|p\|_{\mathcal{H}}^2 = |p(0)|^2 + \|Lp\|_{\mathcal{H}}^2 + \|\Delta^{1/2}Lp\|^2 = |p(0)|^2 + \|Lp\|_{\mathcal{H}}^2 + \|(Lp)(C)v\|^2.
	\]
	Iterating this identity gives
	\begin{equation}\label{norminhbbyee}
		\|p\|_{\mathcal{H}}^2 = \|p\|_{H^2}^2 + \sum_{k \geq 1} \|(L^k p)(C)v\|^2.
	\end{equation}
	By Lemma \ref{graphnormdef}, \eqref{norminhbbyee} implies that $\|p\|_\HH^2=\|p\|_\HG^2$, where $\HG=\HG(E,C,v)$. Since polynomials are dense in both $\HH$ and $\HG$, we conclude that $\HH=\HG$. By Proposition \ref{frknilpappx}, $\HH$ is CNP.
\end{proof}

\begin{remark}
	Theorem \ref{mainthmconcave} provides another proof that every Dirichlet-type space $D(\mu)$ is a CNP space. In general, the space $\HH \in \mathfrak C_{\mathrm{nil}}$ is not a $D(\mu)$ space, see \cite[Theorem 1.2]{LR25}.
\end{remark}

\section{Multipliers and SOT-density of polynomials}\label{sec:Multipliers}

Recall that the multiplier algebra of a RKHS is defined by \eqref{eq:MultiplierAlgebra}. In this section we prove that if $\HH\in\Concave$, and $\psi\in\operatorname{Mult}(\HH)$, then $\psi$ is the SOT-limit of its Cesàro means.

The unitary transformation in Lemma~\ref{hbmodsscnprhB} preserves multiplier symbols and multiplier norms. Hence, we may assume that $\HH=\HH[B]$ for some $B\in\HS_0(\overline{\ran\Delta},\C)$. From Proposition~\ref{prop:Gspaces}, there exist a Hilbert space $E$, a contraction $C$ on $E$, and $v\in E$ such that
\[
\|f\|_{\HH}^2
=
\|f\|_{H^2}^2+\|y(f)\|_{\ell^2(\mathbb N,E)}^2,
\]
where $y(f)$ is the unique solution in $\ell^2(\mathbb N,E)$ to the equations
\begin{equation}\label{norminhbygg}
	y(f)_k-Cy(f)_{k+1}=f_kv,
	\qquad k\ge1.
\end{equation}
In particular, $y\colon\HH\to\ell^2(\mathbb N,E)$ is linear and contractive. Recall that every polynomial belongs to $\HH$. By \eqref{solutionfppc}, for every polynomial $p$,
\[
y(p)_k=\sum_{j\ge0}p_{k+j}C^jv,
\qquad k\ge1.
\]
	
	Fix $n\ge0$, and define
	\[
	a_j=
	\begin{cases}
		1-\dfrac{j}{n+1},&0\le j\le n,\\
		0,&j>n.
	\end{cases}
	\]
	Thus $\sigma_nf=\sum_{j\ge0}a_jf_jz^j$.
	For $1\le k\le n$, telescoping gives
	\begin{align*}
		y(\sigma_nf)_k
		&=\sum_{j=0}^{n-k}a_{k+j}C^jf_{k+j}v\\
		&=\sum_{j=0}^{n-k}a_{k+j}C^j
		\bigl(y(f)_{k+j}-Cy(f)_{k+j+1}\bigr)\\
		&=a_ky(f)_k
		-\frac{1}{n+1}\sum_{j=1}^{n+1-k}C^jy(f)_{k+j}.
	\end{align*}
	For $k>n$, we have $y(\sigma_nf)_k=0$.
	
	For $1\le j\le n$, define
	$B_{j,n}\colon\ell^2(\mathbb N,E)\to\ell^2(\mathbb N,E)$ by
	\[
	(B_{j,n}u)_k=
	\begin{cases}
		C^ju_{k+j},&1\le k\le n+1-j,\\
		0,&k>n+1-j.
	\end{cases}
	\]
	Define $F_n\colon\ell^2(\mathbb N,E)\to\ell^2(\mathbb N,E)$ by
	\[
	(F_nu)_k
	=
	a_ku_k-\frac{1}{n+1}\sum_{j=1}^n(B_{j,n}u)_k.
	\]
	Then $y(\sigma_nf)=F_ny(f)$. Since $\|B_{j,n}\|\le1$ and $0\le a_k\le1$, $\|F_n\|\le2$.
	
	Now let $u\in\ell^2(\mathbb N,E)$ be supported in
	$\{1,\ldots,m\}$. For $n>m$,
	\[
	(F_nu-u)_k=
	\begin{cases}
		-\dfrac{1}{n+1}
		\left(ku_k+\displaystyle\sum_{j=1}^{m-k}C^ju_{k+j}\right),
		&1\le k\le m,\\
		0,&k>m.
	\end{cases}
	\]
	Hence $F_nu\to u$. Since the finitely supported sequences are dense in $\ell^2(\mathbb N,E)$ and $\sup_n\|F_n\|\le2$, we conclude that $F_n\to I$ strongly. Since $\sigma_nf\to f$ in $H^2$, it follows that
	\begin{equation}\label{convergencesignfsc}
		\sigma_nf\to f
		\quad\text{in }\HH
		\qquad(f\in\HH).
	\end{equation}
	
	Since $M_\psi^*k_w=\overline{\psi(w)}k_w$, $\|\psi\|_{H^\infty}\le \|M_\psi\|$. Identify $\ell^2(\mathbb N,E)$ with $H^2(E)$ by sending $(u_k)_{k\ge1}$ to $\sum_{k\ge1}u_kz^{k-1}$, and let $M_\psi^{H^2(E)}$ denote multiplication by $\psi$ on $H^2(E)$. We claim that if $p$ is a polynomial, then
	\begin{equation}\label{ypsipexprpf}
		y(\psi p)
		=
		p(C)y(\psi)+M_\psi^{H^2(E)}y(p),
	\end{equation}
	where $p(C)$ acts coordinatewise on $\ell^2(\mathbb N,E)$. If $\psi$ is also a polynomial, then
	\begin{align*}
		y(\psi p)_k
		&=\sum_{m\ge k}(\psi p)_mC^{m-k}v\\
		&=\sum_{\substack{i,l\ge0\\i+l\ge k}}
		\psi_ip_lC^{i+l-k}v\\
		&=\sum_{i=0}^{k-1}\sum_{l\ge k-i}
		\psi_ip_lC^{l-(k-i)}v
		+\sum_{i\ge k}\sum_{l\ge0}
		\psi_ip_lC^lC^{i-k}v\\
		&=\sum_{i=0}^{k-1}\psi_i y(p)_{k-i}
		+p(C)y(\psi)_k.
	\end{align*}
	For a general $\psi$, take a sequence of polynomials $\psi_n$ converging to $\psi$ in $\HH$. Convergence in $\HH$ implies convergence in \eqref{ypsipexprpf}.
	
	The coordinatewise action of $C$ commutes with $F_n$, and therefore so does the coordinatewise action of $p(C)$.
	Using \eqref{ypsipexprpf} twice, we obtain
	\begin{align*}
		y((\sigma_n\psi)p)
		&=p(C)y(\sigma_n\psi)
		+M_{\sigma_n\psi}^{H^2(E)}y(p)\\
		&=F_n\bigl(p(C)y(\psi)\bigr)
		+M_{\sigma_n\psi}^{H^2(E)}y(p)\\
		&=F_n\left(y(\psi p)-M_\psi^{H^2(E)}y(p)\right)
		+M_{\sigma_n\psi}^{H^2(E)}y(p).
	\end{align*}
	Since $\|F_n\|\le2$, and $\|\sigma_n\psi\|_\infty\le\|\psi\|_\infty\le\|M_\psi\|$, it follows that
	\begin{align*}
		\|y((\sigma_n\psi)p)\|_{\ell^2(\mathbb N,E)}
		&\le
		2\|y(\psi p)-M_\psi^{H^2(E)}y(p)\|_{\ell^2(\mathbb N,E)}
		+\|\sigma_n\psi\|_\infty
		\|y(p)\|_{\ell^2(\mathbb N,E)}\\
		&\le
		2\|y(\psi p)\|_{\ell^2(\mathbb N,E)}
		+3\|M_\psi\|\,
		\|y(p)\|_{\ell^2(\mathbb N,E)}.
	\end{align*}
	Moreover, $\|(\sigma_n\psi)p\|_{H^2}\le\|\psi\|_\infty\|p\|_{H^2}$, so
	\begin{align*}
		\|(\sigma_n\psi)p\|_{\HH}^2
		&=
		\|(\sigma_n\psi)p\|_{H^2}^2
		+\|y((\sigma_n\psi)p)\|_{\ell^2(\mathbb N,E)}^2\\
		&\le
		\|\psi\|_\infty^2\|p\|_{H^2}^2
		+8\|y(\psi p)\|_{\ell^2(\mathbb N,E)}^2
		+18\|M_\psi\|^2\|y(p)\|_{\ell^2(\mathbb N,E)}^2\\
		&\le
		27\|M_\psi\|^2\|p\|_{\HH}^2.
	\end{align*}
Since polynomials are dense in $\HH$, we conclude that $\sup_{n}\|M_{\sigma_n\psi}\|\le\sqrt{27}\|M_\psi\|$. For each polynomial $p$, boundedness of $M_p$ and \eqref{convergencesignfsc} give that $(\sigma_n\psi)p \to \psi p$ in $\HH$, so $M_{\sigma_n\psi}\to M_\psi$ in the strong operator topology.

\section{Applications}\label{sec:Applications}

\subsection{Commutant lifting for concave operators}\label{CLTfcs}

For a general Hilbert space $\HH$, recall that $T \in \mathcal{B}(\mathcal{H})$ is called \emph{analytic} if $\bigcap_{n \geq 0}T^n \HH = \{0\}$. $T$ is called \emph{cyclic}, if there exists $f \in \HH$ such that $\HH=[f]_T$. A result by Richter~\cite{Ri88} ensures that if $T$ is concave, analytic, and $\dim\ker T^*=1$, then $T$ is cyclic.

\begin{theorem}\label{cmtltfcacoh}
	Let $\HH$ be a non-zero Hilbert space, and \(T\in\HB(\HH)\) be cyclic, analytic, and concave. If $\HK\subseteq \HH$ is a closed $T^*$-invariant subspace, and $X\in\HB(\HK)$ commutes with $A=P_\HK T|_\HK$, then there exists $Y\in\HB(\HH)$ such that $YT=TY$, \(P_{\mathcal{K}}Y=XP_{\mathcal{K}}\), and \(\|Y\| = \|X\|\).
\end{theorem}

\begin{proof}
	Since $T$ is expansive, $\ran T$ is closed. Cyclicity implies that $\HH / \ran T$ is at most one-dimensional. Since $T$ is analytic, $\dim \ker T^*=1$, for otherwise $T$ would be invertible. As in the proof of Lemma~\ref{hbmodsscnprhB}, \cite[Theorem~4.6]{LGR} implies that it suffices to consider $T=(M_z,\HH[B])$ for some $\HH[B]\in\ConcaveNorm$. By Theorem~\ref{mainthmconcave}, $\HH[B]$ is CNP. For any polynomial $p$, $\HK$ is invariant under $M_p^*$, and $X$ commutes with the compressed polynomial multiplier $P_\HK M_p|_\HK=p(A)$. By Theorem~\ref{thm:Cesaro}, $\HK$ is invariant under the adjoint of any multiplier, and $X$ commutes with any compressed multiplier. Our conclusion now follows from \cite[Theorem~1.1]{Sh05}.
\end{proof}

\subsection{Singly generated invariant subspaces}

For a second application of Theorem~\ref{thm:Cesaro}, let $\operatorname{Lat}(M_z,\HH)$ be the lattice of closed $M_z$-invariant subspaces of $\HH$. $M_z$-invariant subspaces are natural objects to study. Profound results in this direction include Beurling's theorem for the Hardy space~\cite{Be48}, the Aleman--Richter--Sundberg theorem for the Bergman space~\cite{ARS96}, Richter's and Sundberg's theorems for the Dirichlet space~\cite{Ri88,RS92}. Despite the ongoing and active research on $\HH[B]$ spaces, we know little about the $M_z$-invariant subspaces of $\HH[B]$. In 2019, Aleman and Malman made great progress on this problem. \cite[Theorem~5.11]{AlemanMalman} states that if $\HH=\HH[B]$ is $M_z$-invariant and has finite rank $n$, then every nonzero $\HM \in \operatorname{Lat}(M_z,\HH)$ has index one, meaning $\dim(\HM\ominus z\HM)=1$, and any non-zero $\psi\in\HM\ominus z\HM$ is a cyclic vector for $(M_z,\HM)$, meaning $\HM=[\psi]_{(M_z,\HM)}$. Moreover, there exists a space $\HH[B_1]$ with rank $k\le n$ for which $\HM=\psi\HH[B_1]$. It appears to be an open problem how to describe $\operatorname{Lat}(M_z,\HH[B])$ for general $\HH[B]$ of infinite rank. We extend the Aleman--Malman theorem to infinite rank under the additional hypothesis that $\HH[B]\in\ConcaveNorm$.
\begin{theorem}\label{invariantsubspsigg}
Suppose $\HH[B] \in \ConcaveNorm$. Let $\{0\} \neq \HM \in \operatorname{Lat}(M_z, \HH[B])$. Then \(\mathcal{M} \ominus z\mathcal{M}\) is one-dimensional. If $\psi$ is a unit vector in this space, then
\begin{align*}
\psi \in \operatorname{Mult}(\HH[B]), \qquad \|M_\psi\| = 1, \qquad \HM = [\psi]_{(M_z,\HH[B])} = \overline{\psi\mathcal{H}[B]}.
\end{align*}
Furthermore, there is a Hilbert space \(\mathcal{F}\) and a Schur function \(B_1 \in \mathcal{S}_0(\mathcal{F}, \mathbb{C})\) such that
\[
\mathcal{M} = \psi\mathcal{H}[B_1], \qquad \|\psi g\|_{\mathcal{H}[B]} = \|g\|_{\mathcal{H}[B_1]}, \quad g \in \mathcal{H}[B_1].
\]
The operator \((M_z, \mathcal{H}[B_1])\) is bounded and concave, and the inclusion \(\mathcal{H}[B] \subseteq \mathcal{H}[B_1]\) is contractive.
\end{theorem}
\begin{proof}
Theorem~\ref{mainthmconcave} implies $\HH[B]$ is a CNP space.
Let $\psi\in\operatorname{Mult}(\HH[B])$. By Theorem~\ref{thm:Cesaro}, $\psi \HM \subseteq \HM$. If $\psi(0) = 0$, then $M_{\psi/z} = L M_\psi$. So $\psi/z$ is also a multiplier.

Write \(W = \mathcal{M} \ominus z\mathcal{M}\). If \(W = \{0\}\), then \(\mathcal{M} = z\mathcal{M}\) and \(\mathcal{M} \subseteq \bigcap_{n \geq 0} T^n\mathcal{H}[B] = \{0\}\), contradicting the hypothesis that $\HM\ne \{0\}$. Thus \(W\) contains a unit vector. For any unit vector \(\psi \in W\), \cite[Theorem 1.1]{AHMR17} yields multipliers \(a, b\) such that
\[
\psi = \frac{a}{1-b}, \qquad \|M_a\| \leq 1, \qquad \|M_b\| \leq 1, \qquad b(0) = 0.
\]
Because \(b/z\) is a multiplier and \(\mathcal{M}\) is multiplier invariant, \(b\psi \in z\mathcal{M}\). As \(\psi \perp z\mathcal{M}\), we have
\[
1 \geq \|a\|_{\mathcal{H}[B]}^2 = \|\psi - b\psi\|_{\mathcal{H}[B]}^2 = 1 + \|b\psi\|_{\mathcal{H}[B]}^2.
\]
It follows that \(b\psi = 0\) and \(\psi = a\), so \(\|M_\psi\| \leq 1\). Applying \(M_\psi\) to \(1\) gives \(\|M_\psi\| = 1\).

Choose a nonzero $\psi \in W$ of minimal vanishing order $m$ at zero. For $h \in W$, set $u = \psi/z^m$
and $v = h/z^m$. Both are multipliers and $u(0) \neq 0$. Since $uh = v\psi$,
\[
u(0)h - v(0)\psi = (v - v(0))\psi - (u - u(0))h \in z\mathcal{M} \cap W = \{0\}.
\]
Thus $h$ is a scalar multiple of $\psi$, proving $\dim W = 1$. The restriction $A = M_z|_{\mathcal{M}}$ is analytic and
concave. Then \cite[Theorem 1]{Ri88} gives $\mathcal{M} = [\psi]_{A}$ for any unit vector $\psi \in W$. Polynomial density and boundedness of $M_\psi$ give $[\psi]_{A} = \overline{\psi\mathcal{H}[B]}$.

Now we show $\HM = \psi \HH[B_1]$ for some Schur function $B_1$. Note that $A = M_z|_{\mathcal{M}}$ is expansive and analytic, and \(\ker A^* = W = \mathbb{C}\psi\). Thus by \cite[Theorem 4.6]{LGR}, there exists \(B_1 \in \mathcal{S}_0(\mathcal{F}, \mathbb{C})\) and a unitary \(U \colon \mathcal{H}[B_1] \to \mathcal{M}\) satisfying \(UM_z = AU\), where \(\mathcal{F} = \overline{\operatorname{ran}(A^*A - I)}\). Since \(B_1(0) = 0\), we have \(1 \perp  z\mathcal{H}[B_1]\). Thus \(U\) can be chosen with \(U1 = \psi\). For every polynomial \(p\),
\[
Up = p(A)\psi = \psi p.
\]
Since the functions $\psi p$ span a dense subspace of $\mathcal{M}$, the polynomials are dense in $\mathcal{H}[B_1]$. Thus $\mathcal{M} = \psi \mathcal{H}[B_1]$ with $\|\psi g\|_{\mathcal{H}[B]} = \|g\|_{\mathcal{H}[B_1]}$ whenever $g \in \mathcal{H}[B_1]$. By the concavity of $A$, we obtain that $(M_z, \mathcal{H}[B_1])$ is concave.

For $g \in \mathcal{H}[B]$, we have $\psi g \in \mathcal{M}$ and hence $\psi g = Uh = \psi h$ for some $h \in \mathcal{H}[B_1]$. So $h = g$. Note that
\[
\|g\|_{\mathcal{H}[B_1]} = \|\psi g\|_{\mathcal{H}[B]} \leq \|g\|_{\mathcal{H}[B]}.
\]
Thus \(\mathcal{H}[B] \subseteq \mathcal{H}[B_1]\) and the inclusion is contractive.
\end{proof}

\begin{corollary}
The restriction of a cyclic analytic concave operator to any nonzero closed invariant subspace is cyclic.
\end{corollary}

\subsection{Weak product spaces}
It is a classical theorem that any function in $H^1$ is the product of two functions in $H^2$. For a general RKHS $\mathcal{H}$, consider the weak product space
\begin{equation*}
	\mathcal{H}\odot\mathcal{H}=\left\{f=\sum_{n=1}^\infty \varphi_n g_n \mid \varphi_n,g_n\in\mathcal{H},\  \sum_{n=1}^\infty\|\varphi_n\|_{\mathcal{H}} \|g_n\|_{\mathcal{H}}<\infty\right\}.
\end{equation*}
The weak product space $\mathcal{H}\odot\mathcal{H}$ is a Banach space when equipped with the norm
\begin{equation*}
	\|f\|_{\mathcal{H}\odot\mathcal{H}}=\inf\left\{\sum_{n=1}^\infty\|\varphi_n\|_{\mathcal{H}} \|g_n\|_{\mathcal{H}} \mid
	f=\sum_{n=1}^\infty \varphi_n g_n\right\}.
\end{equation*}
If $\mathcal{H}$ is a CNP space, then $\mathcal{H}\odot\mathcal{H}$ is in many ways the appropriate analogue of $H^1$, see e.g. \cite{AHMR21, AHMR23, ARSW10, ARSW11, CHS21, Har23, JM19, LR15, RS14} and the references therein. If $\mathcal{H}$ is a normalized CNP space, then every $f\in\mathcal{H}\odot\mathcal{H}$ admits a factorization $f=gh$, where $g,h\in\mathcal{H}$ and
\begin{equation*}
	\|f\|_{\mathcal{H}\odot\mathcal{H}}=\|g\|_{\mathcal{H}}\|h\|_{\mathcal{H}}.
\end{equation*}
Moreover, $\operatorname{Mult}(\mathcal{H}\odot\mathcal{H})=\operatorname{Mult}(\mathcal{H})$ and
	\begin{equation*}
		\|\psi\|_{\operatorname{Mult}(\mathcal{H}\odot\mathcal{H})}=\|\psi\|_{\operatorname{Mult}(\mathcal{H})}, \quad \psi\in\operatorname{Mult}(\mathcal{H}),
	\end{equation*}
see \cite[Theorems 1.3, 1.4]{Har23}.

Let $\HH \in \ConcaveNorm$, $X  = \HH \odot \HH$. Every multiplier of $\mathcal{H}$ acts boundedly on $X$ by multiplying one factor in each product. In particular, $M_z$ acts boundedly on $X$. Every $\Lambda\in X^*$ restricts to a bounded functional $h\mapsto \langle h,b\rangle_{\mathcal{H}}$ on $\mathcal{H}$, since $\|1\|_{\mathcal{H}}=1$ and $\|h\|_X\le \|h\|_{\mathcal{H}}$. The associated bounded little Hankel operator $H_b:\mathcal{H}\to \overline{\mathcal{H}}$ is defined by
\[
\langle H_b f,\overline{g}\rangle_{\overline{\mathcal{H}}} = \Lambda(fg) \qquad (f,g\in\mathcal{H}).
\]
Here $\overline{\mathcal{H}}$ is the conjugate Hilbert space, and $|\Lambda(fg)|\le \|\Lambda\|\|f\|_{\mathcal{H}}\|g\|_{\mathcal{H}}$. The following result extends \cite[Theorems 1.1, 1.2, and Corollary 3.3]{LR15}.

\begin{theorem}\label{correinthwpskehkl}
Let $\HH \in \ConcaveNorm$, and $X  = \HH \odot \HH$. Then the maps
\begin{align*}
\operatorname{Lat}(M_z, \mathcal{H}) &\to \operatorname{Lat}(M_z, X),\qquad \mathcal{M} \mapsto \overline{\mathcal{M}}^X, \\
\operatorname{Lat}(M_z, X) &\to \operatorname{Lat}(M_z, \mathcal{H}),\qquad \mathcal{N} \mapsto \mathcal{N} \cap \mathcal{H}
\end{align*}
are inverse lattice isomorphisms. In particular, for every \(f \in \mathcal{H}\),
\[
[f]_{(M_z,X)} \cap \mathcal{H} = [f]_{(M_z,\HH)},
\]
and \(f\) is cyclic in \(\mathcal{H}\) if and only if it is cyclic in \(X\).

For every \(\mathcal{M} \in \operatorname{Lat}(M_z, \mathcal{H})\) there is a bounded little Hankel operator \(H_b\) with \(\mathcal{M} = \ker H_b\). If \(\mathcal{M} \neq \{0\}\), then \(\mathcal{M} \ominus z\mathcal{M}\) is one-dimensional, and for any unit vector \(\psi \in \mathcal{M} \ominus z\mathcal{M}\) one can take \(b = M_z^* \psi\).
\end{theorem}
\begin{proof}
By Theorem~\ref{mainthmconcave}, $\HH$ is CNP. By \cite[Theorem~1.3]{Har23}, every $f\in X$ admits a factorization $f=gh$ with $\|f\|_{\HH\odot\HH}=\|g\|_{\mathcal{H}}\|h\|_{\mathcal{H}}$. By~\cite[Theorem 1.4]{Har23}, $\operatorname{Mult}(X)=\operatorname{Mult}(\mathcal{H})$ in the sense of an isometric isomorphism.

By Theorem~\ref{thm:Cesaro}, polynomials are dense in $\HH$. By the factorization $f=gh$, polynomials are dense in $X$. Moreover, if $\psi\in\operatorname{Mult}(\mathcal{H})$, then $M_{\sigma_n\psi}\to M_\psi$ strongly on $\HH$. For $g,h\in\mathcal{H}$ we have
\[
\|(\sigma_n\psi-\psi)gh\|_X \le \|(\sigma_n\psi-\psi)g\|_{\mathcal{H}}\|h\|_{\mathcal{H}} \to 0.
\]
This implies $M_{\sigma_n\psi}\to M_\psi$ strongly on $X$. Consequently, on both $\mathcal{H}$ and $X$ a closed subspace is $M_z$-invariant if and only if it is multiplier invariant. The weak-product correspondence~\cite[Theorem 1.5]{Har23} thus gives the asserted inverse maps. Both preserve inclusion. The contractive embedding \(\mathcal{H} \subseteq X\) gives \([f]_{(M_z,X)} = \overline{[f]_{(M_z,\HH)}}^X\), so applying the inverse map gives \([f]_{(M_z,X)} \cap \HH = [f]_{(M_z,\HH)}\). Since \(\mathcal{H}\) is dense in \(X\), this also proves the cyclicity equivalence.

Let \(\{0\} \neq \mathcal{M} \in \operatorname{Lat}(M_z, \mathcal{H})\). Let $\psi$ be a unit vector in $\mathcal{M} \ominus z\mathcal{M}$. Theorem \ref{invariantsubspsigg} implies that $\HM = [\psi]_{M_z}$. The weak-product factorization \cite[Theorem 4.1]{Har23}, applied to \(\psi/\|\psi\|_X\), gives \(a, \tilde a \in \operatorname{Mult}(\mathcal{H})\) such that
\[
\psi = \frac{a}{(1-\tilde a)^2}, \qquad \tilde a(0) = 0, \qquad \|M_{\tilde a}\| \leq 1, \qquad \|M_a\| \leq \|\psi\|_X.
\]
Since \(\mathcal{M}\) is multiplier invariant,
\[
a - \psi = z[(\tilde a - 2)(\tilde a/z)\psi] \in z\mathcal{M}.
\]
Orthogonality and the contractive embedding yield
\[
1 \leq (1 + \|a - \psi\|_{\mathcal{H}}^2)^{1/2} = \|a\|_{\mathcal{H}} \leq \|M_a\| \leq \|\psi\|_X \leq 1.
\]
It follows that \(a = \psi\) and \(\|\psi\|_X = 1\).

By the Hahn--Banach theorem, there is a norm-one functional \(\Lambda_\psi \in X^*\) with \(\Lambda_\psi(\psi) = 1\). Its restriction to \(\mathcal{H}\) has norm at most one and takes the value one at the unit vector \(\psi\). Thus
\[
\Lambda_\psi(f) = \langle f, \psi \rangle_{\mathcal{H}}, \qquad f \in \mathcal{H}.
\]
Consequently, \(H_\psi\) is bounded. Since multiplication by \(z\) is bounded on \(X\), the functional \(h \mapsto \Lambda_\psi(zh)\) defines a bounded Hankel operator with symbol \(b = T^*\psi\), and \(H_b = H_\psi T\), where $T = (M_z, \HH)$.

For \(f \in \mathcal{H}\) and any polynomial \(p\),
\[
\langle H_b f, \overline{p} \rangle_{\overline{\mathcal{H}}} = \langle zpf, \psi \rangle_{\mathcal{H}}.
\]
If \(f \in \mathcal{M}\), then \(zpf \in z\mathcal{M}\), so \(\mathcal{M} \subseteq \ker H_b\) by polynomial density. Put \(\mathcal{N} = \ker H_b\). The displayed identity also shows that \(\mathcal{N}\) is \(T\)-invariant. For \(f \in \mathcal{N}\), taking \(p = 1\) gives \(\langle Tf, \psi \rangle_{\mathcal{H}} = 0\). Hence \(\psi \in \mathcal{N} \ominus z\mathcal{N}\). Applying Theorem \ref{invariantsubspsigg} to \(\mathcal{N}\) gives
\[
\mathcal{N} = [\psi]_{(M_z,\mathcal{H})} = \mathcal{M}.
\]

For $\mathcal{M}=\{0\}$, choose distinct $r_j\in(0,1/2)$ tending to zero and positive $c_j$ such that $\sum_j c_j K^\HH(r_j,r_j)<\infty$. The estimate $|h(r_j)|\le K^\HH(r_j,r_j)\|h\|_X$ shows that $\Lambda(h)=\sum_j c_j h(r_j)$ defines a bounded functional on $X$. Let $b$ be its symbol. Then for polynomial $p$,
\[
\langle H_b f, \overline{p} \rangle_{\overline{\mathcal{H}}} = \Lambda(fp) = \sum_j c_j f(r_j) p(r_j).
\]
If $f\in\ker H_b$, then $\sum_j c_j f(r_j)p(r_j)=0$. Taking $p$ to be the Taylor polynomials of $z\mapsto \overline{f(\bar z)}$ and using uniform convergence on $[0,1/2]$ gives $\sum_j c_j |f(r_j)|^2=0$. Hence $f$ vanishes at every $r_j$ and is identically zero.
\end{proof}

\begin{remark}
(i) Jury constructed a CNP space on \(\mathbb{D}\) with bounded \(M_z\) and a closed \(M_z\)-invariant subspace that fails to be multiplier invariant, see \cite[Theorem~3 and Section~3]{Jur05}. For a Hankel operator $H_b$, \(\ker H_b\) is multiplier invariant, as is any intersection of such kernels. It follows that Jury's \(M_z\)-invariant subspace cannot equal the kernel of a bounded little Hankel operator, or even an intersection of such kernels. Consequently, the CNP property and boundedness of \(M_z\) alone do not imply the Hankel conclusion of Theorem \ref{correinthwpskehkl} for all closed \(M_z\)-invariant subspaces.

(ii) For a normalized CNP space $\HH$, Hartz \cite[Theorem 4.4]{Har23} proved that every multiplier invariant subspace $\HM$ is the countable intersection of kernels of Hankel operators. In general, $\HM$ may not equal the kernel of a single Hankel operator, see the example of Richter--Sunkes~\cite[Example 4.3]{RSu16} in the Drury--Arveson space setting.
\end{remark}

By applying the above theorem, we obtain the following approximation result for the $M_z$-invariant subspaces, cf.~\cite[Corollary~4.5]{LR15}.
\begin{corollary}
Let $\HH \in \ConcaveNorm$. Let \(\mathcal{M}_n\) and \(\mathcal{M}\) be nonzero closed \(M_z\)-invariant subspaces of \(\mathcal{H}\). Let \(\psi_n \in \mathcal{M}_n \ominus z\mathcal{M}_n\) and \(\psi \in \mathcal{M} \ominus z\mathcal{M}\) be unit vectors. If \(\psi_n(z) \to \psi(z)\) for every \(z \in \mathbb{D}\), then \(P_{\mathcal{M}_n} \to P_{\mathcal{M}}\) in the strong operator topology.
\end{corollary}
\begin{proof}
Since $\psi_n(z)\to\psi(z)$ and $\|\psi_n\|_\HH=\|\psi\|_\HH=1$, $\psi_n\to\psi$ weakly by compactness. Consequently, $\|\psi_n-\psi\|_\HH^2=2-2\mathrm{Re}\langle\psi_n,\psi\rangle_\HH\to 0$. \cite[Proposition~4.3]{LR15} implies that $\HM=[\psi]_{M_z}\subseteq \underline{\lim}\, \HM_n$, i.e., for every $f\in\HM$ there exists $f_n\in\HM_n$ such that $f_n\to f$ in $\HH$.

By Theorem~\ref{correinthwpskehkl}, $\HM_n=\ker H_{b_n}$, where $b_n=M_z^*\psi_n$. A similar identification holds for $\HM$. Clearly $b_n\to b$ in $\HH$. By \cite[Lemma~2.2]{LR15}, $\HM=[b]_*^\perp$, where $[b]_*$ is the smallest closed subspace containing $b$ which is invariant under all $M_\phi^*$, $\phi\in\operatorname{Mult}(\HH)$. It follows from Theorem~\ref{thm:Cesaro} that $[b]_*=[b]_{M_z^*}$. Consequently, $\HM^\perp\subseteq \underline{\lim}\, \HM_n^\perp$. The SOT convergence of $P_{\HM_n}$ follows from \cite[Lemma~4.2]{LR15}.
\end{proof}

\section*{Acknowledgement}
We, the authors, thank Alexandru Aleman, Michael Hartz and Stefan Richter for many stimulating conversations. We acknowledge the use of ChatGPT~6 in the preparation of this paper. Collaboration between the authors began in 2021, and we have previously published \cite{LR25}. We have had many initial ideas for the present paper, including using the de~Branges--Rovnyak spaces, the construction in Section~\ref{sec:Comparison}, and the rank-one extension scheme introduced by Stefan Richter and used in \cite[Section 2.2]{GL}. Through discussions with ChatGPT, we proved Lemma~\ref{rankoneextens}, which led us to prove Theorem~\ref{mainthmconcave}. ChatGPT was also used for the proofreading and literature search stages of the paper. All mathematical statements, proofs, and references were independently checked, revised, and finalized by the authors, who take full responsibility for the content of this paper.

S. Luo was supported by National Natural Science Foundation of China (No. 12271149), Natural Science Foundation of Hunan Province (No. 2024JJ2008).

E. Rydhe was supported by the Swedish Research Council, grant 2022-04307.

\end{document}